\documentclass[11pt]{article}
\usepackage{amsfonts}
\usepackage{latexsym}
\usepackage{graphicx}

\def\filleftmap{\mathord\leftarrow \mkern-6mu
	\cleaders\hbox{$\mkern-2mu \mathord- \mkern-2mu$}\hfill
	\mkern-6mu \mathord-}

\newtheorem{thm}{Theorem}[section]
\newtheorem{defn}[thm]{Definition}
\newtheorem{lemma}[thm]{Lemma}

\newtheorem{cor}[thm]{Corollary}
\newtheorem{prop}[thm]{Proposition}

\newtheorem{example}[thm]{Example}
\newtheorem{remark}[thm]{Remark}
\def\<{{\langle}}
\def\>{{\rangle}}
\def\Ast{{\textstyle \ast}}

\def\Irr{{\rm Irr}{\scriptscriptstyle \,}}
\def\Min{{\rm Min}{\scriptscriptstyle \,}}
\def\invlim#1{{{\displaystyle \lim_{\longleftarrow}}\;\!_{#1}\,}}
\def\Invlim#1{{\displaystyle \lim_{\stackrel{\longleftarrow}{{#1}}}\,}}
\def\from{{\,\leftarrow\,}}
\def\hdot{{^{\textstyle \cdot}}}
\def\ldot{{.{\scriptstyle \,}}}
\def\hhom{{{\underline{\rm Hom}}{\scriptscriptstyle \,}}}
\def\eext{{{\underline{\rm Ext}}{\scriptscriptstyle \,}}}
\def\ttor{{{\underline{\rm Tor}}{\scriptstyle \;}}}
\def\aa{{{\rm \bf a}}}
\def\bb{{{\rm \bf b}}}
\def\cc{{{\rm \bf c}}}
\def\ee{{{\rm \bf e}}}
\def\ww{{{\rm \bf w}}}
\def\mm{{\mathfrak m}}
\def\pp{{\mathfrak p}}

\def\HH{{\widetilde H}{\scriptstyle \!}{}}

\def\FF{{\mathbb F}}

\def\NN{{\mathbb N}}
\def\RR{{\mathbb R}}
\def\ZZ{{\mathbb Z}}

\begin{document}

\centerline{{\large \bf Alexander Duality for Monomial Ideals}
		{\large \bf and Their Resolutions}}

\vskip 3mm
\centerline{Ezra Miller}
\vskip 1mm
\begin{abstract}
\noindent
Alexander duality has, in the past, made its way into commutative algebra
through Stanley-Reisner rings of simplicial complexes.  This has the
disadvantage that one is limited to squarefree monomial ideals.  The
notion of Alexander duality is generalized here to arbitrary monomial
ideals.  It is shown how this duality is naturally expressed by Bass
numbers, in their relations to the Betti numbers of a monomial ideal and
its Alexander dual.
Relative cohomological constructions on cellular complexes are shown to
relate cellular free resolutions of a monomial ideal to free resolutions
of its Alexander dual ideal.  As an application, a new canonical
resolution for monomial ideals is constructed.
\end{abstract}

{}

\vskip 1mm
\centerline{{\bf Introduction}}
\vskip .5mm Alexander duality in its most basic form is a relation
between the homology of a simplicial complex $\Gamma$ and the cohomology
of another simplicial complex $\Gamma^\vee$, called the \emph{dual} of
$\Gamma$.  Recently there has been much interest (\cite{TH}, \cite{ER},
\cite{HRW}, \cite{BCP}) in the consequences of this relation when applied
to the monomial ideals which are the Stanley-Reisner ideals $I_\Gamma$
and $I_{\Gamma^\vee}$ for the given simplicial complex and its Alexander
dual.  This has the limitation that Stanley-Reisner ideals are always
squarefree.  The first aim of this paper is to define Alexander duality
for arbitrary monomial ideals and then generalize some of the relations
between $I_\Gamma$ and $I_{\Gamma^\vee}$.  A second goal is to
demonstrate that Bass numbers are
the proper vessels for the translation of Alexander duality into
commutative algebra.  The final 
goal is to reveal the connections between Alexander duality and the
recent work on cellular resolutions.

There are two ``minimal'' ways of describing an arbitrary monomial ideal:
via the minimal generators or via the (unique) irredundant irreducible
decomposition.  Given a monomial ideal $I$, Definition~\ref{defn:alexdual}
describes a method for producing another monomial ideal $I^\vee$ whose
minimal generators correspond to the irredundant irreducible components
of $I$.  Miraculously, this is enough to guarantee that the minimal
generators of $I$ correspond to the irreducible components of $I^\vee$.
It is particularly easy to verify that this reversal of roles takes place
for the squarefree ideals $I = I_\Gamma$ and $I^\vee = I_{\Gamma^\vee}$
above (Proposition~\ref{prop:X-dual}).  A connection with linkage and
canonical modules is described in Theorem~\ref{thm:linkage}.

One can also deal with Alexander duality as a combinatorial phenomenon,
thinking of $\Gamma$ as an order ideal in the lattice of subsets of
$\{1,\ldots,n\}$.  The Alexander dual $\Gamma^\vee$ is then given by the
complement of the order ideal, which gives an order ideal in the opposite
lattice.  For squarefree monomial ideals all is well since the only
monomials we care about are represented precisely by the lattice of
subsets of $\{1,\ldots,n\}$.  For general monomial ideals we instead
consider the larger lattice $\ZZ^n$, by which we mean the poset with its
natural partial order $\preceq$.  Then a monomial ideal $I$ can be
regarded as a dual order ideal in $\ZZ^n$, and $I^\vee$ is constructed
(roughly) from the complementary set of lattice points, which is an order
ideal---see Definition~\ref{defn:[a]}.  It is Theorem~\ref{thm:[a]=a}
which proves the equivalence of the two definitions.

Bass numbers first assert themselves in Section~\ref{section:bass-betti}.
Their relations to Betti numbers for monomial modules
(Corollary~\ref{cor:lattice} and Theorem~\ref{thm:gor}) are derived as
consequences of graded local duality and Alexander duality (in its avatar
as lattice duality in $\ZZ^n$).  The Bass-Betti relations are then
massaged to equate the localized Bass numbers of $I$
(Definition~\ref{defn:bass}) with the Betti numbers of $I^\vee$ in the
first of the two central results of this paper,
Theorem~\ref{thm:bass-betti}.  Theorem~\ref{thm:[a]=a} is then recovered
as a special case of this main result, which also finds an application to
inequalities between the Betti numbers of dual ideals
(Theorem~\ref{thm:inequalities}) generalizing those for squarefree ideals
in \cite{BCP}.

The extension of Alexander duality to resolutions is accomplished in
Sections~\ref{section:cellular} and~\ref{section:limits}.  A new
canonical and geometric resolution, the {\it cohull resolution} is
constructed in Definition~\ref{defn:cohull}.  It should be thought of as
Alexander dual to the {\it hull resolution} of \cite{BS} (which is
similarly canonical and geometric).  Roughly speaking, the cohull
resolution is constructed from the irreducible components instead of the
minimal generators.  The cohull resolution owes its existence to the
second central result of the paper, Theorem~\ref{thm:relcocell}, which is
a more general result on duality for cellular resolutions.  Its proof,
which is resolutely algebraic, is the content of
Section~\ref{section:limits}.  The idea is to deform an ideal into its
dual step by step via Definition~\ref{defn:fb} and keep track of the
deformations on cellular resolutions (Theorem~\ref{thm:limit}).  The
final step, taken in Theorem~\ref{thm:coresolution}, is to check the
effect of the deformations on the homology of the resolutions.

\vskip 2mm
\noindent
{\bf Acknowledgements.} The author would like to express his thanks to
Dave Bayer, 
David Eisenbud, 
Serkan Ho\c sten, 
Sorin Popescu, 
Stefan Schmidt, 
Frank Sottile, 
Bernd Sturmfels, 
and Kohji Yanagawa 
for their helpful comments and discussions.

\begin{section}{Definitions and basic properties}

\label{section:definitions}

For notation, let $S$ be the $\ZZ^n$-graded $k$-algebra
$k[x_1,\ldots,x_n] \subseteq T := S[x_1^{-1}, \ldots, x_n^{-1}]$, where
$k$ is a field and $n \geq 2$.  If $A \subseteq T$ is any subset, $\<a
\mid a \in A\>$ will denote the $S$-submodule generated by the elements
in $A$, and it may also be regarded as an ideal if $A \subseteq S$.  The
maximal $\ZZ^n$-graded ideal $\<x_1,\ldots,x_n\>$ of $S$ will be denoted
by $\mm$.  Each (Laurent) monomial in $T$ is specified uniquely by a
single vector $\aa = (a_1,\ldots,a_n) = \sum_i a_i \ee_i \in \ZZ^n$,
while each irreducible monomial ideal is specified uniquely by a vector
$\bb = (b_1,\ldots,b_n) \in \NN^n$, so the notation
$$
  x^\aa = x_1^{a_1} \cdots x_n^{a_n} \qquad {\rm and} \qquad
  \mm^\bb = \<x_i^{b_i} \mid b_i \geq 1\>
$$ 
will be used to highlight the similarity.  The $\ZZ^n$-graded prime
ideals, which are precisely the monomial prime ideals, are indexed by
faces of the $(n-1)$-simplex $\Delta := 2^{\{1,\ldots,n\}}$ with vertices
$1,\ldots,n$.  Identifying a face $F \in \Delta$ with its characteristic
vector in $\ZZ^n$, the monomial prime corresponding to $F$ may be written
with the above notation as $\mm^F$.  Note, in particular, that $\mm^\bb$
need not be an artinian ideal, just as $x^\aa$ need not have full
support.  In fact, $\mm^\bb$ is $\mm^{\sqrt\bb}$-primary, where $\sqrt\bb
\in \Delta$ is the face representing the support of $\bb$; that is,
$\sqrt\bb$ has $i^{\,\rm th}$ coordinate 1 if $b_i \geq 1$ and 0
otherwise.  With this notation, taking radicals can be expressed as
$\sqrt{\mm^\bb} = \mm^{\sqrt\bb}$.

All modules $N$ and homomorphisms of such will be $\ZZ^n$-graded, so that
$N = \bigoplus_{\aa \in \ZZ^{\scriptstyle n}} N_\aa$.  In addition, any
module that is isomorphic to a submodule of $T$ \emph{$\!\!\!$ as a
$\ZZ^n$-graded module} will, if it is convenient, be freely identified
with that submodule of $T$.  For instance, the principal ideal generated
by $x_1 \cdots x_n$ can be identified with the module $S[-{\bf 1}]$,
where ${\bf 1} = (1,\ldots,1) \in \ZZ^n$ and $N[\aa]_\bb = N_{\aa+\bb}$
for $\aa, \bb \in \ZZ^n$.  In this paper, ideals will all be proper
monomial ideals, and the symbol $I$ will always denote such an ideal.
The vector $\aa_I$ will denote the exponent on the least common multiple
of the minimal generators of $I$.

Before making the definition of Alexander dual ideal, the next few
results make sure that the exponents used to define the set $\Irr(I)$ of
irredundant irreducible components of $I$ are $\preceq \aa_I$.  For the
next two results, let $\Lambda$ denote the set of irreducible ideals
containing $I$.

\begin{lemma}
If $\mm^\bb \in \Irr(I)$ then $\mm^\bb$ is minimal (under inclusion) in
$\Lambda$.
\end{lemma}
{\it Proof:\ } Suppose $\mm^\bb \neq \mm^\cc$ and that $\mm^\bb \supseteq
\mm^\cc \in \Lambda$.  If now $I = \mm^\bb \cap I'$ for some ideal $I'$
then also $I = \mm^\cc \cap I'$, whence $\mm^\bb \not\in \Irr(I)$.
\hfill
$\Box$
\begin{prop}
If $\mm^\bb \in \Irr(I)$ then for each $i \in \sqrt\bb$ there is a
minimal generator $x^\cc$ of $I$ with $b_i = c_i$.
\end{prop}
{\it Proof:\ } Suppose $\mm^\bb \in \Irr(I)$ but the conclusion does not
hold.  Then given any minimal generator $x^\cc$ of $I$, either $b_{i'}
\leq c_{i'}$ for some $i \neq i' \in \sqrt\bb$, or else $b_i < c_i$.  In
either case, $x^\cc \in \mm^{\bb + \ee_i}$, where $\ee_i$ is the
$i^{\,\rm th}$ unit vector in $\ZZ^n$.  Then $\mm^{\bb + \ee_i} \supseteq
I$, contradicting the minimality of $\mm^\bb$ in $\Lambda$.
\hfill
$\Box$
\begin{cor} \label{cor:Irr}
For any $\mm^\bb \in \Irr(I)$ we have $\bb \preceq \aa_I$. \hfill $\Box$
\end{cor}

The following notation will be very convenient in the definition and
handling of Alexander duality.  For any vector $\aa \in \ZZ^n$ and any
face $F \in \Delta$, let $\aa \cdot F$ denote the restriction of $\aa$ to
$F$:
$$
  (\aa \cdot F)_i\ =\ \left\{\begin{array}{ll}
				a_i & {\rm if}\ i \in F \cr
				0     & {\rm otherwise}
			\end{array}\right..
$$
This operation may also be thought of as the coordinatewise product of
$\aa$ and $F$.  If, in addition, ${\bf 0} \preceq \bb \preceq \aa$,
define $\bb^\aa$ to be the vector whose $i^{\,\rm th}$ coordinate is $a_i
+ 1 - b_i$ if $b_i \geq 1$ and $0$ otherwise; more compactly,
$$
  \bb^\aa\ =\ (\aa + {\bf 1} - \bb) \cdot \sqrt\bb
	 \ =\ (\aa + {\bf 1}) \cdot \sqrt\bb - \bb\,,
$$
where $\sqrt\bb$ is the support of $\bb$, as above.  The next result is a
first indication of the utility of $\bb^\aa$ when applied to irreducible
ideals $\mm^\bb$.

\begin{prop} \label{prop:irr_order}
If \,${\bf 0} \preceq \bb,\cc \preceq \aa$\, then \,$\mm^\bb \supseteq
\mm^\cc$ if and only if \,$\bb^\aa \succeq \cc^\aa$.
\end{prop}
{\it Proof:\ } The condition $\mm^\bb \supseteq \mm^\cc$ is equivalent to
the combination of (i) $\sqrt\bb \succeq \sqrt\cc$ and (ii) $\bb \cdot
\sqrt\cc \preceq \cc$.  Now consider the inequalities in the following
chain:
$$
	\bb^\aa\ \ =\ \ (\aa + {\bf 1} - \bb) \cdot \sqrt\bb
	\ \ \succeq\ \  (\aa + {\bf 1} - \bb) \cdot \sqrt\cc
	\ \ \succeq\ \  (\aa + {\bf 1} - \cc) \cdot \sqrt\cc
	\ \ =\ \ \cc^\aa\,.
$$
The left inequality is equivalent to (i) since $\aa + {\bf 1} - \bb$ has
full support, and the right inequality is equivalent to (ii) since $\cc
\cdot \sqrt\cc = \cc$.  It remains only to show that $\bb^\aa \succeq
\cc^\aa$ implies both inequalities, and this can be checked
coordinatewise.  If $c_i = 0$, then both inequalities become trivial; if
$c_i > 0$ then $b_i > 0$, and the left inequality becomes an equality
while the right inequality becomes $(\bb^\aa)_i = a_i + 1 - b_i \geq a_i
+ 1 - c_i = (\cc^\aa)_i$\,.  \hfill $\Box$
\vskip 2mm

Corollary~\ref{cor:Irr} clears the way for the main definition of this
paper:

\begin{defn}[Alexander duality] \label{defn:alexdual}
Given an ideal $I$ and $\aa \succeq \aa_I$, the \emph{Alexander dual
ideal $I^\aa$ with respect to $\aa$} is defined by
$$
	I^\aa \quad = \quad \<x^{\bb^{\scriptstyle \aa}} \mid \mm^\bb
	\in \Irr(I)\>.
$$
For the special case when $\aa = \aa_I$, let $I^\vee = I^{\aa_I}$.
\end{defn}

\begin{remark} \label{rk:alexdual} \rm
(i) We will never have occasion to take an Alexander dual of the ideal
$\mm$, so $\mm^\aa$ will retain its original definition.
\vskip 1mm

\noindent
(ii) The dual $I^\aa$ with respect to any $\aa \succeq \aa_I$ depends
only on $\aa \cdot \sqrt{\aa_I}$.  This is because $\bb$ and $\aa \cdot
\sqrt{\bb}$ determine $\bb^\aa$, and $\aa \cdot \sqrt\bb = (\aa \cdot
\sqrt{\aa_I}) \cdot \sqrt\bb$ for all of the relevant $\bb$ by
Corollary~\ref{cor:Irr}.  In particular, $I^\vee = I^{\bf 1}$ if $I$ is
squarefree.
\vskip 1mm

\noindent
(iii) $I^\vee$ is not gotten by taking the depolarization of the
Alexander dual of the polarization of $I$ (see \cite{SV}, Chapter~II for
polarization).  For instance, when $I = \<x^2, xy, y^2\>$, the
polarization is $I_{polar} = \<x_1x_2, x_1y_1, y_1y_2\>$, whose canonical
Alexander dual is $I_{polar}^\vee = \<x_1y_1, x_1y_2, x_2y_1\>$.
Removing the subscripts on $x$ and $y$ then yields the principal ideal
$\<xy\>$, whereas $I^\vee = \<xy^2, x^2y\>$.
\end{remark}

\begin{prop} \label{prop:min}
The set of generators for $I^\aa$ given by the definition is minimal.
More generally, suppose $\aa \succeq \aa_I$ and $\Lambda$ is a collection
of integer vectors $\preceq \aa$ such that $I = \bigcap_{\bb \in
\Lambda}\,\mm^\bb$.  Then $I^\aa = \<x^{\bb^{\scriptstyle \aa}} \mid \bb
\in \Lambda\>$, and the intersection determined by $\Lambda$ is
irredundant if and only if the set of generators for $I^\aa$ is minimal.
\end{prop}
{\it Proof:\ } This follows from Corollary~\ref{cor:Irr} and
Proposition~\ref{prop:irr_order}.
\hfill
$\Box$
\begin{example}\rm \label{ex:alexdual}
\begin{figure}[htb] 
{\centering
	\includegraphics{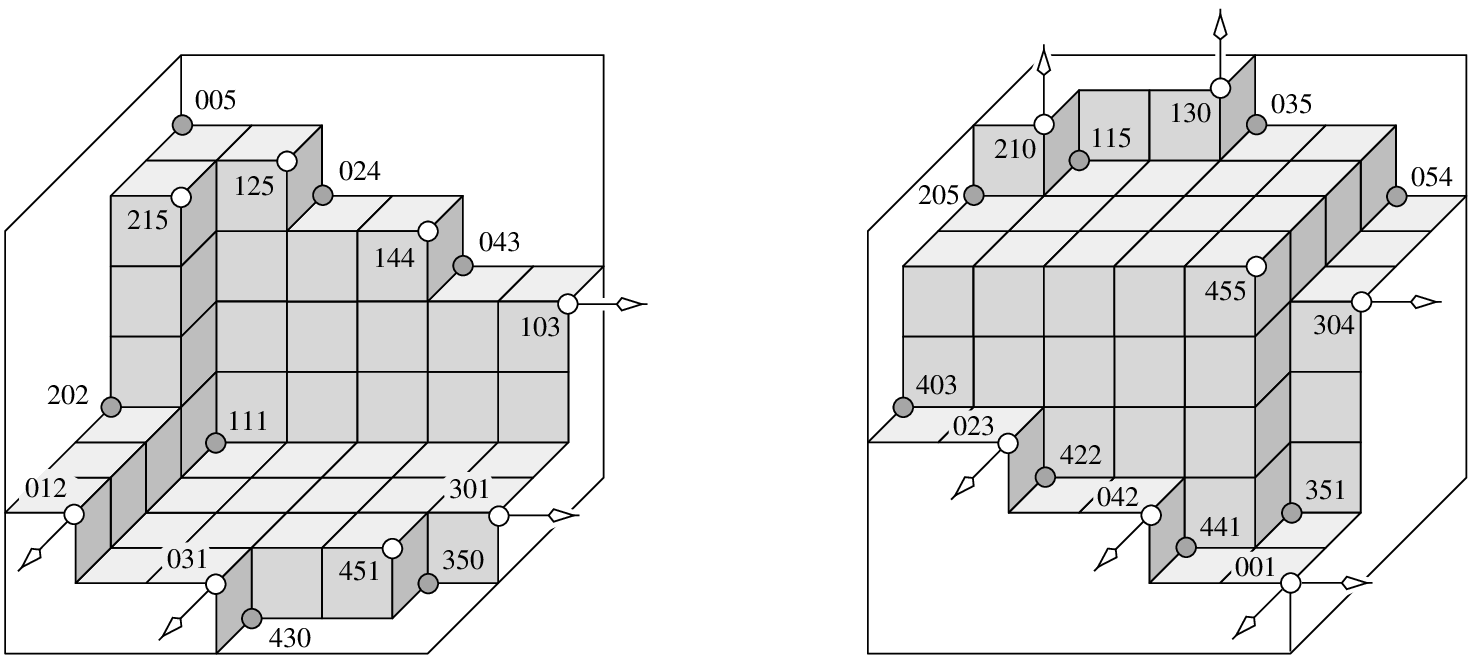}\\
	\hskip .25in $I$ \hskip 3.2in $I^\vee$ \hskip .55in \mbox{}\\
}
$$
\begin{array}{rcl}
I
& = &	\<z^5, x^2z^2, x^4y^3, x^3y^5, y^4z^3, y^2z^4, xyz\>
\cr
& = &	\<x^2,y,z^5\> \cap \<y,z^2\> \cap \<y^3,z\> \cap \<x^4,y^5,z\>
	\cap \<x^3,z\> \cap \<x,z^3\> \cap \<x,y^4,z^4\> \cap
	\<x,y^2,z^5\>
\cr
\end{array}
$$
\hskip 1cm $\aa := \aa_I = (4,5,5)$
$$
\begin{array}{rcl}
I^\vee
& = &

	\<z\> \cap \<x^3,z^4\> \cap \<x,y^3\> \cap \<x^2,y\> \cap
	\<y^2,z^3\> \cap \<y^4,z^2\> \cap \<x^4,y^5,z^5\> \phantom{so\
	that\ so\ that\ }

\cr
& = &	\<x^3y^5z, y^5z^4, y^3z^5, xyz^5, x^2z^5, x^4z^3, x^4y^2z^2,
	x^4y^4z\>.
\cr
\end{array}
$$
\caption{\label{fig:stairs} \baselineskip 0pt 
The truncated staircase diagrams, minimal generators, and irredundant
irreducible components for $I$ and $I^\vee$.  Black lattice points are
generators, and white lattice points indicate irreducible components.
The numbers are to be interpreted as vectors, e.g.\ 205 = (2,0,5).  The
arrows attached to a white lattice point indicate the directions in which
the component continues to infinity; it should be noted that a white
point has a zero in some coordinate precisely when it has an arrow
pointing in the corresponding direction.}
\end{figure}
Let $n = 3$, so that $S = k[x,y,z]$.  Figure~\ref{fig:stairs} lists the
minimal generators and irredundant irreducible components of an ideal $I
\subseteq S$ and its dual $I^\vee$ with respect to $\aa_I$.  The
(truncated) ``staircase diagrams'' representing the monomials not in
these ideals are also rendered in Figure~\ref{fig:stairs}.  In fact, the
staircase diagram for $I^\vee$ is gotten by literally turning the
staircase diagram for $I$ upside-down (the reader is encouraged to try
this).  Notice that the support of a minimal generator of $I$ is equal to
the support of the corresponding irreducible component of $I^\vee$.
\hfill
$\Box$
\end{example}
\begin{example} \label{ex:permut} \rm
Let $\Sigma_n$ denote the symmetric group on $\{1,\ldots,n\}$ and $\cc =
(1,2,\ldots,n) \in \NN^n$.  The ideal $I = \<x^{\sigma(\cc)} \mid \sigma
\in \Sigma_n\>$ is the {\it permutahedron ideal} determined by $\cc$,
introduced in \cite{BS}, Example~1.9.  The results of
Example~\ref{ex:permut3} below imply that the canonical Alexander dual is
the {\it tree ideal}, which is generated by $2^n - 1$ monomials: $I^\vee
= \<(x^F)^{n-|F|+1} \mid \emptyset \neq F \in \Delta\>$.  For instance,
when $n = 3$,
\begin{eqnarray*}
I\; & = & \<xy^2z^3,xy^3z^2,x^2yz^3,x^2y^3z,x^3yz^2,x^3y^2z\> \cr
I^\vee & = & \<xyz,x^2y^2,x^2z^2,y^2z^2,x^3,y^3,z^3\>.
\end{eqnarray*}
The tree ideal is so named becuase it has the same number $(n+1)^{n-1}$
of standard monomials (monomials not in the ideal) as there are trees on
$n+1$ labelled vertices.  The minimal free resolution of $I^\vee$ is
obtained in Example~\ref{ex:permut3}, below. \hfill $\Box$
\end{example}

Recall that for a simplicial complex $\Gamma \subseteq \Delta$ the
\emph{Stanley-Reisner ideal} $I_\Gamma$ of $\Gamma$ is defined by the
nonfaces of $\Gamma$:
$$
\begin{array}{rcl}
\phantom{\ \<x^F \mid F \not\in \Gamma\>,}
\phantom{\Gamma^\vee\ }
\phantom{I_\Gamma\ }
	I_\Gamma
	& = &
	\<x^F \mid F \not\in \Gamma\>,
\phantom{\ \{\overline{F} \in \Delta \mid F \not \in \Gamma\}.}
\phantom{\bigcap_{F \in \Gamma} \mm^{\overline F}}
\end{array}
$$
and the \emph{Alexander dual simplicial complex} $\Gamma^\vee$ consists
of the complements of the nonfaces of $\Gamma$:
$$
\begin{array}{rcl}
\phantom{\ \<x^F \mid F \not\in \Gamma\>,}
\phantom{I_\Gamma\ }
\phantom{I_\Gamma\ }
	\Gamma^\vee
	& = &
	\{F \in \Delta \mid {\overline F} \not \in \Gamma\},
\phantom{\ \<x^F \mid F \not\in \Gamma\>,}
\phantom{\bigcap_{F \in \Gamma} \mm^{\overline F}}
\end{array}
$$
where ${\overline F} = \{1,\ldots,n\} \setminus F$.  Recall also that
$I_\Gamma$ may be equivalently described as
$$
\begin{array}{rcl}
\phantom{\ \<x^F \mid F \not\in \Gamma\>,}
\phantom{\Gamma^\vee\ }
\phantom{I_\Gamma\:}
	I_\Gamma
	& = &
	{\displaystyle \bigcap_{{\overline F} \in \Gamma} \mm^F},
\phantom{\ \<x^F \mid F \not\in \Gamma\>,}
\phantom{\ \{\overline{F} \in \Delta \mid F \not \in \Gamma\}.}
\end{array}
$$
since $\mm^F \supseteq I\ \Leftrightarrow\ F$ has at least one vertex in
each nonface of $\Gamma\ \Leftrightarrow\ {\overline F}$ is missing at
least one vertex from each nonface of $\Gamma\ \Leftrightarrow\
{\overline F}$ is a face of $\Gamma$.  Applying
Definition~\ref{defn:alexdual} to the latter characterization of
$I_\Gamma$ yields:

\begin{prop} \label{prop:X-dual}
For a simplicial complex $\Gamma \subseteq \Delta$ we have $I_\Gamma^\vee
= I_{\Gamma^\vee}$.
\end{prop}
{\it Proof:\ } Observe that $\bb^{\bf 1} = \bb$ if $\bb \in \{0,1\}^n$,
and use Proposition~\ref{prop:min} along with
Remark~\ref{rk:alexdual}(ii).  We get $I_\Gamma^\vee\ =\ \<x^F \mid
{\overline F} \in \Gamma\>\ =\ \<x^F \mid F \not\in \Gamma^\vee\>\ =\
I_{\Gamma^\vee}$.
\hfill
$\Box$

Thus, as promised, Definition~\ref{defn:alexdual} generalizes to
arbitrary monomial ideals the definition of Alexander duality for
squarefree monomial ideals.  The connection with the squarefree case is
never lost, however, because the general definition does the same thing
to the zero-set of $I$ as the squarefree definition does:
\begin{prop}
Taking Alexander duals commutes with taking radicals: $\sqrt{I^\vee} =
\sqrt{I}^{\,\vee}$.
\end{prop}
{\it Proof:\ } Since ${\bf 0} \preceq \bb \preceq \aa_I$ whenever
$\mm^\bb \in \Irr(I)$, the equality $\sqrt{\bb} = \sqrt{\bb^{\aa_I}}\,$
follows from the definitions.  Thus,
$$
\begin{array}{rcl}
\sqrt{I^\vee} &
	= &
	\<x^{\sqrt \bb} \mid \mm^\bb \in \Irr(I)\> \cr
& = &	\<x^F \mid \mm^F {\rm \ is\ minimal\ among\ primes\ containing\ }
	I\>\cr
& = &	\sqrt{I}^{\,\vee}, \cr
\end{array}
$$
the last equality using again the facts mentioned in the first line of
the proof of Proposition~\ref{prop:X-dual}.
\hfill
$\Box$
\vskip 2mm

The notion of Alexander duality sheds light on the interconnections
between some of the developments in \cite{BPS}, \cite{BS}, and \cite{Stu}
concerning cellular resolutions and (co)generic monomial ideals.  To
begin with, consider the following condition on a set of vectors $\{\bb^j
= (b^j_1,\ldots,b^j_n) \in \NN^n\}_{j=1}^r$:
$$
b_i^j \geq 1\ \Rightarrow\ b^j_i \neq b^{j'}_i\ {\rm for\ all\ } j' \neq j.
$$
A \emph{generic} ideal, as defined in \cite{BPS}, is an ideal whose
minimal generators have exponent vectors satisfying the above condition;
similarly, a \emph{cogeneric} ideal, as defined in \cite{Stu}, is an
ideal whose irredundant irreducible components have exponent vectors
satisfying the above condition.  Using Definition~\ref{defn:alexdual} the
following is immediate (for any $\aa \succeq \aa_I$).
\begin{prop} \label{prop:gen-cogen}
$I^\aa$ is generic if and only if $I$ is cogeneric.
\hfill
$\Box$
\end{prop}
\begin{example}\rm
The ideal $I$ in Example~\ref{ex:alexdual} is generic, while $I^\vee$ is
cogeneric.
\hfill
$\Box$
\end{example}
The connections between the minimal resolutions of such ideals and
cellular resolutions will be explored in Section~\ref{section:cellular}.

Recall that the {\it Castelnuovo-Mumford regularity} and {\it initial
degree} of a $\ZZ$-graded $S$-module $L$ defined respectively by
$$
  {\rm reg}(L) := \max\{j \in \ZZ \mid {\rm Tor}_i(L,k)_{i + j} \neq 0\}
  \qquad {\rm and} \qquad
  {\rm indeg}(L) := \min\{j \in \ZZ \mid L_j \neq 0\},
$$
where $L_j$ is the $j^{\rm th}$ component of $L$.  The question was
raised in \cite{HRW}, Question~10 whether there is a duality for possibly
nonradical monomial ideals with the ``amazing properties''
$$
\begin{array}{rl}
\bullet &
	\rm{reg}(I) - \rm{indeg}(I) = \dim (S/I^\vee) -
	\rm{depth}(S/I^\vee)
\cr
\bullet &
	I \hbox{ is componentwise linear if and only if } S/I^\vee
	\hbox{ is sequentially Cohen-Macaulay}
\end{array}
$$
obeyed by Alexander duals in the squarefree case.  Here, $I$ is
considered in its $\ZZ$-grading.  Having defined a duality operation in
this paper, some comments are obviously warranted.

First of all, it is unrealistic to expect the first property to extend to
the arbitrary (nonradical) case since the right-hand side of the equation
is bounded while the left-hand side is not, in general.  For instance, if
$d \in \NN$ then $\ {\rm reg}(\mm^{d \cdot {\bf 1}}) - {\rm indeg}(\mm^{d
\cdot {\bf 1}}) = n(d-1) - d\ $ while $(\mm^{d \cdot {\bf 1}})^\vee =
\<x_1 \cdots x_n\>$ is Cohen-Macaulay.  Nevertheless, there may be some
class of ideals which behaves nicely under some kind of duality, not
necessarily as defined here.  As to whether or not such a class of ideals
exists for the Alexander duality as defined here, such an investigation
has not yet been made.

Unfortunately, the second property also fails for $I$ and $I^\aa$, for
somewhat trivial reasons: almost every ideal has an artinian Alexander
dual.  Specifically, if $I$ is arbitrary and $x = x_1 \cdots x_n$, then
$S/(xI)^\aa$ is artinian (for any $\aa \succeq \aa_I$), and hence
Cohen-Macaulay.  But the minimal free resolution of $xI$ is just the
shift by ${\bf 1}$ of the minimal resolution of $I$.  Thus every minimal
resolution, be it componentwise linear or not, appears as the resolution
of an ideal whose dual is a Cohen-Macaulay ideal; i.e.\ $S/I^\aa$
Cohen-Macaualy $\not \Rightarrow$ $I$ componentwise linear.

One might still hope that the implication ``$I$ has a linear resolution
$\Rightarrow S/I^\aa$ is sequentially Cohen-Macaulay'' would hold, but
even this fails, as the example below shows.  The fundamental problem
with the nonsquarefree case is that the $\ZZ$-degree of an element is not
determined by the support of its $\ZZ^n$-graded degree, as it is with
squarefree monomials.  Thus an ideal might have a linear resolution while
its generators have support sets of varying sizes, wreaking havoc with
the equidimensionality required for the Cohen-Macaulayness of the dual.
Even so, it would be very interesting to know what is the property
Alexander dual to ``sequentially Cohen-Macaulay''; perhaps this property
could relax the requirements of componentwise linearity in a nice way.

\begin{example}\rm
Let $I' = \<ab,bc,cd\> \subseteq S = k[a,b,c,d]$ be the ideal of the
``stick twisted cubic'' simplicial complex spanned by the edges $\{b,d\},
\{b,c\}$, and $\{a,c\}$.  It is readily checked that $I'$ has a linear
resolution: indeed, $(I')^\vee$ is the ideal of another stick twisted
cubic, which is Cohen-Macaulay because the stick twisted cubic is
connected and has dimension 1, so \cite{ER}, Theorem~3 applies.  Let
$$
\begin{array}{l}
  I = \mm I' = \<a^2b,abc,acd,ab^2,b^2c,bcd,abc,bc^2,c^2d,abd,bcd,cd^2\>
\cr
  I^\vee = \<b^2d^2,b^2c^2,a^2c^2,abc^2d^2,a^2bcd^2,a^2b^2cd\>
\end{array}
$$
with $\aa_I = (2,2,2,2)$.  Then $I$ has a linear resolution by
\cite{HRW}, Lemma~1, and we show that $S/I^\vee$ is not sequentially
Cohen-Macaulay.

Recall that for a module $N$ to be sequentially Cohen-Macaulay, we
require that there exist a filtration $0 = N_0 \subset N_1 \subset
\cdots \subset N_r = N$ such that $N_i/N_{i-1}$ is Cohen-Macaulay for
all $i \leq r$ and $\dim(N_{i+1}/N_{i}) > \dim(N_{i}/N_{i-1})$ for all
$i < r$.  It follows from the equidimensionality of $N/N_{r-1}$ and
the strict reduction of dimension in successive quotients that
$N_{r-1}$ is the top dimensional piece of $N$; i.e.\ $N_{r-1}$ is the
intersection of all primary components (of $0$ in $N$) which have
dimension $\dim(N)$.  Thus it suffices to check that $S/I^\vee_{\rm
top}$ is not Cohen-Macaulay, where $I^\vee_{\rm top} = \<b^2d^2,b^2cd,
abcd,b^2c^2,abc^2,a^2c^2\>$ is the intersection of all primary
components of $I^\vee$ which have dimension $2 = \dim (S/I^\vee)$.
\hfill $\Box$
\end{example}

\end{section}
\begin{section}{Alternate characterizations of the Alexander dual ideal}%

\label{section:alternate}

\noindent
Definition~\ref{defn:alexdual} is quite satisfactory for the consequences
just derived from it, but it can sometimes be inconvenient to work with.
For instance, it is not obvious from the definition that $(I^\aa)^\aa =
I$, which is fundamental---see Corollary~\ref{cor:Iaa=I}.  For this and
other applications, we set out now to find other characterizations of the
Alexander dual ideal in Theorem~\ref{thm:linkage} and in
Definition~\ref{defn:[a]} with Theorem~\ref{thm:[a]=a}.  Along the way,
an algebraic analogue of combinatorial lattice duality in $\ZZ^n$ is
defined in Defintion~\ref{defn:T-dual}.

First, a result relating Alexander duality to linkage (see \cite{Vas},
Appendix~A.9 for a brief introduction to linkage and references, or
\cite{SV} for more details):
\begin{thm} \label{thm:linkage}
If $\,\aa \succeq \aa_I\,$ then $(\mm^{\aa+{\bf 1}} : I^\aa) = I +
\mm^{\aa + {\bf 1}}$.
\end{thm}
{\it Proof:\ } Let $\Min(I^\aa)$ denote the exponents on the minimal
generators of $I^\aa$.  Then $(\mm^{\aa + {\bf 1}} : I^\aa) =
\bigcap_{\bb \in \Min(I^\aa)} (\mm^{\aa + {\bf 1}} : x^\bb)$.  But $x^\cc
\cdot x^\bb \in \mm^{\aa + {\bf 1}} \ \Leftrightarrow\ \bb + \cc
\not\preceq \aa \ \Leftrightarrow\ \cc \not\preceq \aa - \bb \
\Leftrightarrow\ x^\cc \in \mm^{\aa + {\bf 1} - \bb}$.  Thus, taking all
intersections over $\bb \in \Min(I^\aa)$,
$$
  \bigcap\,(\mm^{\aa + {\bf 1}} : x^\bb)
\ =\ 
  \bigcap\,\mm^{\aa + {\bf 1} - \bb}
\ =\
  \bigcap\,\Bigl( \mm^{\bb^{\scriptstyle \aa}} + \mm^{\aa + {\bf 1}}
  \Bigr)
\ =\
  \Bigl( \bigcap\,\mm^{\bb^{\scriptstyle \aa}} \Bigr) + \mm^{\aa + {\bf
  1}}
\ =\
  I + \mm^{\aa + {\bf 1}}
$$
since $(\bb^\aa)^\aa = \bb$ for all $\bb \preceq \aa$.
\hfill
$\Box$
\vskip 2mm 
\noindent
\begin{remark} \rm
Using Corollary~\ref{cor:Iaa=I}, below, this theorem provides a useful
way to compute the Alexander dual ideal, given a set of generators.
Indeed, the generators for $I^\aa$ are simply those generators of
$(\mm^{\aa + {\bf 1}} : I)$ whose exponents are $\preceq \aa$.  Using
Definition~\ref{defn:alexdual} (and Corollary~\ref{cor:Iaa=I} again),
this can also be construed as a method for computing irreducible
components of $I$ given a generating set for $I$, or vice versa.
\end{remark}

Denoting the $\ZZ^n$-graded Hom functor by $\hhom$, the next duality that
comes into play is the {\it $k$-dual} $N^\wedge := \hhom_k(N,k)$, which
is a $\ZZ^n$-graded $S$-module with the grading $(N^\wedge)_\cc = {\rm
Hom_{\,}}_k (N_{-\cc},k)$.  It is a simple but very important observation
that $T^\wedge \cong T$ as $\ZZ^n$-graded modules.  This can be
exploited: let $M \subseteq T$ be a submodule (the $\ZZ^n$-graded
submodules of $T$ are precisely the {\it monomial modules} of \cite{BS}).
Taking the $k$-dual of the surjection $T \to T/M$ yields an injection
$(T/M)^\wedge \to T^\wedge \cong T$.  This makes $(T/M)^\wedge$ into a
submodule of $T$ which we call the {\it $T$-dual of $M$} and denote by
$M^T$.  If one thinks of the module $M$ as a set of lattice points in
$\ZZ^n$, then $M^T$ can be thought of as the negatives of the lattice
points in the complement of $M$; i.e.\ we can make the equivalent
\begin{defn} \label{defn:T-dual}
The $T$-{\rm dual} $M^T$ of a monomial module $M \subseteq T$ is defined
by $x^{-\bb} \in M^T \Leftrightarrow x^\bb \not\in M$.
\end{defn}
The equivalence with the earlier formulation can be seen simply by
examining which $\ZZ^n$-graded pieces of $M$ and $M^T$ have dimension $1$
over $k$ and which have dimension $0$.  Observe the striking similarity
of Definition~\ref{defn:T-dual} with definition of the dual simplicial
complex: $\overline{F} \in \Gamma^\vee \Leftrightarrow F \not\in \Gamma$.
Here are some properties of the $T$-dual which will be used later
(possibly without explicit reference).  Note the similarity of (i)--(iii)
to the laws governing complements, unions, and intersections.

\begin{prop} \label{prop:T-dual}
Let $M$ and $N$ be submodules of $T$.  Then
$$
\begin{array}{llcll}
(i)	& (M^T)^T = M &
\qquad&	(v)	& T/M^T = M^\wedge \cr
(ii)	& M \subseteq N \ \Leftrightarrow \ N^T \subseteq M^T &
\qquad&	(vi)	& (N/M)^\wedge = M^T/N^T\ {\rm if\ } M \subseteq N \cr
(iii)	& (M+N)^{\,T} = M^T \cap N^T &
\qquad&	(vii)	& (N/\,N\!\cap\!M)^\wedge = M^T /\,M^T\!\cap\!N^T \cr
(iv)	& M[\aa]^T = M^T[-\aa] \cr
\end{array}
$$
\end{prop}
{\it Proof:\ } Statements (i)--(iv) follow from
Definition~\ref{defn:T-dual}, and (v) follows either from the definition
and (i) or as a special case of (vi).  To prove (vi) observe that $N/M =
\ker(T/M \to T/N)$ so that $(N/M)^\wedge = {\rm coker}((T/N)^\wedge \to
(T/M)^\wedge)$ and use the definition of $T$-dual.  Finally, (vii) is
just (vi) and (iii) applied to $(N+M)/M = N/\,N\!\cap\!M$.
\hfill
$\Box$
\begin{defn}
Given a monomial ideal $I \subseteq S$ define the {\rm \v Cech hull of
$I$} in $T$:
$$
	\widetilde{I}\ :=\ \< x^\bb \mid \bb \in \ZZ^n\ {\rm and}\
	x^{\bb^+}\! \in I \>\,,
$$
where $\bb^+ \in \NN^n$ is, as usual, the join (componentwise maximum) of
$\bb$ and ${\bf 0}$ in the order lattice $\ZZ^n$.
\end{defn}

\begin{prop} \label{prop:tilde}
Taking the \v Cech hull commutes with finite intersections and sums.
Furthermore,
$$
\begin{array}{ll}
(i)&
	\hbox{$\widetilde I$ is the largest monomial submodule of $T$ whose
	intersection with $S$ is equal to $I$.} \cr
(ii)&
	\hbox{$\widetilde I$ can be generated by (possibly infinitely
	many) monomials in $T$ of degree $\, \preceq \aa_I$.} \cr
(iii)&
	\hbox{$\widetilde{I}^{\;T}$ is generated in degrees $\preceq {\bf
	0}$.} \cr
\end{array}
$$
\end{prop}
{\it Proof:\ } The first statement follows from (i) and the definitions.

\vskip 1mm \noindent 
(i) It is clear from the definition that $\widetilde{I}$ contains $I$;
and if $x^\bb \in \widetilde{I} \cap S$ then $\bb^+ = \bb$ whence $x^\bb
\in I$.  Thus $\widetilde{I} \cap S = I$.  On the other hand, if $M$ is a
monomial submodule of $T$ satisfying $M \cap S = I$ and $x^\bb \in M$,
then $x^{\bb^-} \!\!\cdot x^\bb = x^{\bb^+} \in M \cap S = I$, where
$\bb^- := \bb^+ - \bb$.  Thus $M \subseteq \widetilde{I}$.

\vskip 1mm \noindent
(ii) If $x^\bb \in \widetilde{I}$ then $\cc \preceq \bb^+$ for some
minimal generator $x^\cc$ of $I$, whence $x^{\cc - \bb^-}$ is in
$\widetilde{I}$, divides $x^\bb$, and has exponent $\,\preceq \aa_I$.

\vskip 1mm \noindent
(iii) The following statement is precisely the $T$-dual to statement (i):
$\widetilde{I}^{\;T}$ is the smallest submodule whose sum with
$\widetilde{\mm}$ is equal to $I^{\;T}\!$.  As $\widetilde{\mm}$ already
contains all degrees $\not\preceq {\bf 0}$, minimality of
$\widetilde{I}^{\;T}$ implies that it is generated in degrees $\preceq
{\bf 0}$.
\hfill
$\Box$

\begin{example}\rm \label{ex:tilde}
(i)\ Recall that for $F \in \Delta$, the complement $\{1,\ldots,n\}
\setminus F$ is denoted by $\overline F$.  Using this, the localization
$S[x^{-{\overline F}\,}]$ is achieved by inverting the variables $x_i$
for $i \not\in F$.  Now let $\bb \succ {\bf 0}$ and $F = \sqrt\bb$.  Then
$$
  \Bigl(\widetilde{\mm^\bb}\Bigr)^T\ 
	=\ \Bigl(S[x^{-{\overline F}\,}]\Bigr)[\bb-F]\,.
$$
To see this, first observe that if $\cc \in \NN^n$ then $x^\cc \not\in
\mm^\bb\ \Leftrightarrow\ \cc \cdot F \preceq \bb-F$.  Therefore, if $\cc
\in \ZZ^n$ then $x^\cc \not\in \widetilde{\mm^\bb}\ \Leftrightarrow\
\cc^+ \preceq \bb-F\ \Leftrightarrow\ \cc \cdot F \preceq \bb-F$.  This
last condition is equivalent to $-\cc \cdot F \succeq F-\bb$, and this
occurs if and only if $x^{-\cc} \in \Bigl(S[x^{-{\overline F}\,}]\Bigr)
[\bb-F]$.

\noindent
(ii)\ For a special case, it follows that when $\bb \succeq {\bf 0}$, 
$\,\widetilde{\mm^{\bb + {\bf 1}}}\,=\,S[\bb]^T$.
\hfill
$\Box$
\end{example}
\begin{remark} \rm
Example~\ref{ex:tilde}(ii) is the reason for the name {\it \v Cech hull}:
when $\bb = {\bf 0}$, we find that $\widetilde{\mm}$ is the kernel of the
last map in the \v Cech complex on $x_1,\ldots,x_n$.
\end{remark}

\begin{defn} \label{defn:[a]}
For any monomial ideal $I$ and $\aa \succeq \aa_I$, define
$$
	I^{[\aa]}\ :=\ \widetilde{I}^{\;T}[-\aa] \cap S.
$$
\end{defn}
\begin{example} \rm \label{ex:zz2}
Figure~2 is a schematic diagram depicting the transformation in stages
from $I$ to $I^{[\aa_I]}$.  The black and white dots shift by ${\bf 1}$
from the penultimate stage to the last; they are left in place (with
respect to the dark black dot and the dark dotted lines) for the rest of
the transformation.  This shift is the reason for the ${\bf 1}$ in the
definition of $\bb^\aa$, and it occurs because the flip-flop from
$\widetilde{I}$ to $\widetilde{I}^{\:T}$ leaves a space of ${\bf 1}$.
The crux of this whole theory is that the ``boundaries'' of
$\widetilde{I}$ and $\widetilde{I}^{\:T}$ have the same shape, but
reversed, thus switching the roles of the black and white dots.  This
schematic may be helpful in parsing the proof of Theorem~\ref{thm:[a]=a},
below.  \hfill $\Box$
\end{example}
\begin{figure}[!t]
\vbox{\centering
\includegraphics{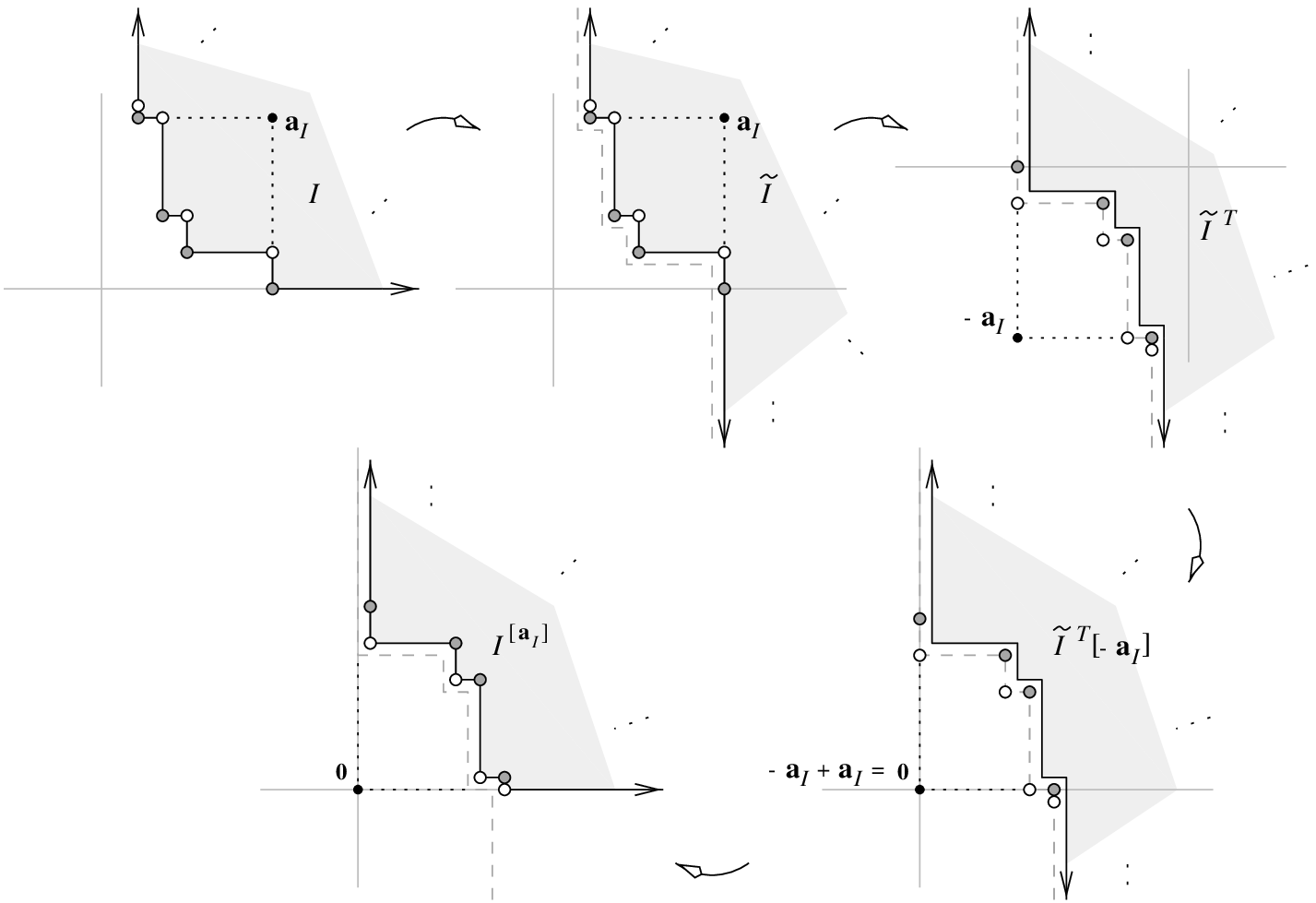}\\
Figure~2}
\end{figure}
\begin{lemma} \label{lemma:Iatilde}
$(I^{[\aa]})\!\widetilde{\phantom{ii}} \! = \widetilde{I}^{\;T}[-\aa]$.
\end{lemma}
{\it Proof:\ } Let $M = \widetilde{I}^{\;T}[-\aa]$.  By
Proposition~\ref{prop:tilde}(i), $(M \cap S) \!  \widetilde
{\phantom{ii}}\!  \supseteq M$ since their intersections with $S$ are
equal by definition.  Thus $((M \cap S)\!\widetilde {\phantom{ii}})^{\,T}
\subseteq M^T$, with equality in degrees $\,\preceq {\bf 0}$.  But $M^T =
\widetilde{I}[\aa]$ is generated in negative degrees by
Proposition~\ref{prop:tilde}(ii), so that in fact $((M \cap
S)\!\widetilde{\phantom{ii}})^{\,T} = M^T$.  Taking $T$-duals of this
equality gives the desired result.
\hfill
$\Box$
\vskip 2mm

\noindent
The upshot is that $I$ may be reconstructed from $I^{[\aa]}$ via the same
construction which produces $I^{[\aa]}$ from $I$:

\begin{prop} \label{prop:Iaa=I}
$\aa_{I^{[\aa]}} \preceq \aa$ and $I = (I^{[\aa]})^{[\aa]}$.
\end{prop}
{\it Proof:\ } By Proposition~\ref{prop:tilde}(iii)
$\widetilde{I}^{\;T}[-\aa]$ is generated in degrees $\,\preceq \aa$, so
Lemma~\ref{lemma:Iatilde} implies that the same holds for
$(I^{[\aa]})\!\widetilde{\phantom{ii}}$.  It is trivial to check that if
any monomial module $M \subseteq T$ is generated in degrees $\preceq \aa$
then so is $M \cap S$, because $\aa \succeq {\bf 0}$.  Thus
$\aa_{I^{[\aa]}} \preceq \aa$, and, in particular, $(I^{[\aa]})^{[\aa]}$
is well-defined.  Now

\vskip -.5\baselineskip
\vspace \baselineskip
\noindent
\parbox{.95\linewidth}{\centering
$
\begin{array}{r@{\ \ =\ \ }llr}
(I^{[\aa]})^{[\aa]} 
&
	((I^{[\aa]})\!\widetilde{\phantom{ii}})^{T}[-\aa] \cap S
&
	\hbox{by definition}
\cr
&
	(\widetilde{I}^{\;T}[-\aa])^T[-\aa] \cap S
&
	\hbox{by the previous lemma} \cr
&
	\widetilde{I} \cap S
&
	\hbox{by Proposition~\ref{prop:T-dual}(iv) and~(i)}
\cr
&
	I.
&
&
\cr
\end{array}
$}
\newlength{\Mylen}
\settodepth{\Mylen}{
\parbox{.95\linewidth}{\centering
$
\begin{array}{r@{\ \ =\ \ }llr}
(I^{[\aa]})^{[\aa]} 
&
	((I^{[\aa]})\!\widetilde{\phantom{ii}})^{T}[-\aa] \cap S
&
	\hbox{by definition}
\cr
&
	(\widetilde{I}^{\;T}[-\aa])^T[-\aa] \cap S
&
	\hbox{by the previous lemma} \cr
&
	\widetilde{I} \cap S
&
	\hbox{by Proposition~\ref{prop:T-dual}(iv) and~(i)}
\cr
&
	I.
&
&
\cr
\end{array}
$}
}
\addtolength{\Mylen}{-1.3mm}
\hfill \raisebox{-\Mylen}{$\Box$}
\vskip 2mm

The real cause for introducing $I^{[\aa]}$ is the next result, which may
not be so unexpected at this point.  It would seem that
Theorem~\ref{thm:[a]=a} makes the notation $I^{[\aa]}$ superfluous, and
it does; nevertheless, the notation will be retained for emphasis, to
indicate that Sections~\ref{section:bass-betti}
and~\ref{section:localize} (and, in particular,
Theorem~\ref{thm:bass-betti}) are logically independent from
Theorem~\ref{thm:[a]=a}.
\begin{thm} \label{thm:[a]=a}
$I^\aa = I^{[\aa]}$.
\end{thm}
{\it Proof:\ } To simplify notation, declare that $\bb \in \Irr(I)$ if
$\mm^\bb \in \Irr(I)$.  For each $\bb \in \Irr(I)$, let $S^\bb =
S[x^{\,\mbox{--}\,\overline{\sqrt{\bb}}\,}]$ be the localization of $S$
at the prime $\mm^{\sqrt{\bb}}$.  Then for each $\bb \in \Irr(I)$ and any
$\cc \in \NN^n$ we have the following two facts:

\noindent
(i) $S^\bb[-\cc] \cong S^\bb[-\cc \cdot \sqrt{\bb}^{\,}]$ since
multiplication by $x^{\cc \cdot \overline{\sqrt{\bb}}}$ is a
$\ZZ^n$-graded automorphism of $S^\bb[-\cc]$.

\vskip 1mm \noindent
(ii) $S \cap S^\bb[-\cc \cdot \sqrt{\bb}^{\,}] = S[-\cc \cdot
\sqrt{\bb}^{\,}]$.  Indeed, this is equivalent to $\Bigl(\<x^{\cc \cdot
\sqrt{\bb}\,}\> \cdot S^\bb\Bigr) \cap S = \<x^{\cc \cdot
\sqrt{\bb}\,}\>$, which holds because $\<x^{\cc \cdot \sqrt{\bb}\,}\>
\subseteq S$ is saturated with respect to
$\<x^{\overline{\sqrt{\bb}}\,}\>$; i.e.\ $\Bigl(\<x^{\cc \cdot
\sqrt{\bb}\,}\> : x^{\overline{\sqrt{\bb}}\,}\Bigr) = \<x^{\cc \cdot
\sqrt{\bb}\,}\>$.

\vskip 1mm
Creating $I^{[\aa]}$ from $I$ in stages yields
$$
\begin{array}{rr@{\quad =\quad }ll}
&
\widetilde{I}
&	\displaystyle \bigcap\,\widetilde{\mm^{\bb}}
&	\hbox{by Proposition~\ref{prop:tilde}}\phantom{\Bigl(\Bigr)}
\cr
\Rightarrow
&
\widetilde{I}^{\;T}
&	\displaystyle \sum \Bigl(\widetilde{\mm^\bb}\Bigr)^T
&	\hbox{by Proposition~\ref{prop:T-dual}(iii)}
\cr
&
&	\displaystyle \sum S^\bb[\bb - \sqrt\bb]
&	\hbox{by Example~\ref{ex:tilde}(i)}\phantom{\Bigl(\Bigr)^T}
\cr
\Rightarrow
&
\widetilde{I}^{\;T}[-\aa]
&	\displaystyle \sum S^\bb[-\bb^\aa]
&	\hbox{by (i) above, with}\ \cc = \aa + \sqrt\bb - \bb
		\phantom{\Bigl(\Bigr)^T}
\cr
\Rightarrow
&
S \cap \widetilde{I}^{\;T}[-\aa]
&	\displaystyle \sum S[-\bb^\aa]
&	\hbox{by (ii) above, with}\ \cc = \bb^\aa
		\phantom{\Bigl(\Bigr)^T}
\cr
\end{array}
$$
where the intersection and all of the summations are taken over all $\bb
\in \Irr(I)$.  The last sum above is equal to $I^\aa$ since each summand
$S[-\bb^\aa]$ is just a principal ideal $\<x^{\bb^{\scriptstyle \aa}}\>$.
\hfill
$\Box$

\begin{cor} \label{cor:Iaa=I}
$(I^\aa)^\aa = I$.  Furthermore, $(\bb^\aa)^\aa = \bb$, so that

\noindent
\parbox{.95\linewidth}{\centering
$ \displaystyle
  I^\aa \quad = \quad \bigcap\,\{\mm^{\bb^{\scriptstyle \aa}} \mid 
			x^\bb\ \hbox{is a minimal generator of}\ I\}\,.
$
}
\settodepth{\Mylen}{
\parbox{.95\linewidth}{\centering
$ \displaystyle
  I^\aa \quad = \quad \bigcap\,\{\mm^{\bb^{\scriptstyle \aa}} \mid 
			x^\bb\ \hbox{is a minimal generator of}\ I\}\,.
$
}}
\addtolength{\Mylen}{-1.3mm}
\hfill \raisebox{-\Mylen}{$\Box$}
\end{cor}
\begin{remark} \rm
In general, one has $(I^\vee)^\vee \neq I$.  However, in the special case
when $I = \sqrt{I}$, it will always happen that $(I^\vee)^\vee = I$.
This follows from an application of Corollary~\ref{cor:Iaa=I} to
Remark~\ref{rk:alexdual}(ii).  The difference $\aa_I - \aa_{I^\vee}$
measures the extent to which $(I^\vee)^\vee \neq I$ fails, in the sense
that $(I^\vee)^\vee = I[\aa_I - \aa_{I^\vee}]\cap S$.  However
$((I^\vee)^\vee)^\vee = I^\vee$, so that an ideal which is already an
Alexander dual is maximal in some sense.  It is unclear what the
invariant $\aa_I - \aa_{I^\vee}$ means, in general, although the passage
from $I$ to $(I^\vee)^\vee$ can sometimes be thought of as a
``tightening'' that may resolve some amount of nonminimality in the hull
resolution of \cite{BS}---see Example~\ref{ex:tighten}.  See also
Remark~\ref{rk:relcocell}(ii) below for another occurrence of the
invariant $\aa_I - \aa_{I^\vee}$.
\end{remark}

The reader interested in cellular resolutions may wish to skip directly
to Section~\ref{section:cellular}, whose only logical dependence on
Sections~\ref{section:bass-betti} and~\ref{section:localize} is
Proposition~\ref{prop:ext}.

\end{section}
\begin{section}{Bass numbers versus Betti numbers}

\label{section:bass-betti}

\noindent
Algebraically, Alexander duality is best expressed in terms of relations
between Betti and Bass numbers (Definition~\ref{defn:bass-betti}), as
evidenced by this section and the next.  The principle behind this is
that the $T$-duality of Section~\ref{section:alternate}, which can be
thought of as lattice duality in $\ZZ^n$, can also be interpreted
(Corollary~\ref{cor:lattice}) as a manifestation of the same process that
interchanges flat and injective modules (in the appropriate category).
In Theorem~\ref{thm:gor} this results in equalities between Bass and
Betti numbers of $I$.  Though perhaps not so interesting a statement in
its own right, Proposition~\ref{prop:ext} is the workhorse for the
remainder of the paper---it is the {\it reason} everything else is true.
It encapsulates simultaneously the relations between all of the dualities
that enter into this paper: $k$-\ and $T$-duality, Alexander duality,
linkage, local duality, and Matlis duality.

\begin{defn} \label{defn:bass-betti}
The derived functors of the $\ZZ^n$-graded functor $\hhom$ will be called
$\eext$, and the left derived functor of $\otimes$, which is also
$\ZZ^n$-graded, will be called $\ttor$.  For a module $N$ define
\begin{eqnarray*}
  \mu_{i,\bb}(N)   & = & \dim_k \Bigl(\eext^i_S(k,N)_\bb \Bigr) \cr
  \beta_{i,\bb}(N) & = & \dim_k \Bigl(\ttor_i^S(k,N)_\bb \Bigr),
\end{eqnarray*}
the {\rm $i^{\rm \,th}$ Bass and Betti numbers of $N$ in degree $\bb$}.
\end{defn}
\begin{remark} \label{rk:derived} \rm
(i) In order to compute these derived functors in the category $\cal M$
of $\ZZ^n$-graded $S$-modules (see Proposition~\ref{prop:calculate}), we
need to know that $\cal M$ has enough injective and projective modules,
just as in the nongraded case.  Of course, there are always free modules,
so this takes care of the projectives; for injectives one can easily
modify the proof of \cite{BH}, Theorem~3.6.2 to fit the $\ZZ^n$-graded
case.

\vskip 1mm \noindent
(ii) If $M$ is finitely generated then $\eext^{\textstyle \cdot}(M,-) =
{\rm Ext\,}^{\textstyle \cdot}(M,-)$.  In particular, summing the Betti
or Bass numbers over all $\bb$ (or all $\bb$ with fixed $\ZZ$-degree)
gives the same result as computing directly in the nongraded (or
$\ZZ$-graded) case.
\end{remark}

In what follows, we will need the notion of a {\it flat resolution in
$\cal M$}.  This is defined exactly like a free resolution, except that
the resolving modules are required to be flat instead of free, where {\it
flat} means acyclic for $\ttor$ \cite{Wei}, Section~2.4.  Recall that
free and flat are equivalent for finitely generated $S$-modules; this is
a simple consequence of the grading and Nakayama's lemma.  However,
non-finitely generated flat modules, such as localizations of $S$, may
fail to be free, or even projective.
\begin{prop} \label{prop:calculate}
(i) $\eext^{\textstyle \cdot}(M,N)$ can be calculated as the homology of
the complexes obtained either by applying $\hhom(-,N)$ to a projective
resolution of $M$ in $\cal M$ or by applying $\hhom(M,-)$ to an injective
resolution of $N$ in $\cal M$. \\
(ii) $\ttor_{\textstyle \cdot}(M,N)$ can be calculated as the homology of
the complexes obtained by either tensoring with $N$ a flat resolution of
$M$ in $\cal M$ or by tensoring with $M$ a flat resolution of $N$ in
$\cal M$.
\end{prop}
{\it Proof:\ } (i) Remark~\ref{rk:derived}(i) above provides enough
injectives to use \cite{Wei}, Definition~2.5.1, Example~2.5.3, and
Exercise~2.7.4.

\vskip 1mm \noindent
(ii) \cite{Wei}, Theorem~2.7.2 and Exercise~2.4.3.  
\hfill
$\Box$

\begin{lemma}
$N^\wedge = \hhom_S(N,S^\wedge)$.
\end{lemma}
{\it Proof:\ } \cite{BH}, Proposition~3.6.16(c), whose proof holds just
as easily in the $\ZZ^n$-graded case.
\hfill 
$\Box$
\vskip 2mm

The next theorem is the starting point for the comparison of Betti and
Bass numbers.  Its corollary, which carries out the lattice
complementation, is fundamental to the rest of the results in this
section.

\begin{prop} \label{prop:complementation}
For any module $N$, 
$
	\mu_{i,\bb}(N)\ =\ \beta_{i,-\bb}(N^\wedge).
$
\end{prop}
{\it Proof:\ } A module $J$ is injective if and only if $J^\wedge$ is
flat, because 
\begin{equation} \label{eqn:dash}
  \hhom(-\,,\,J) \quad = \quad \hhom(-\,,\,\hhom(J^\wedge,S^\wedge)\,)
  \quad = \quad \hhom(- \otimes J^\wedge\,,\,S^\wedge).
\end{equation}
Indeed, the first term being an exact functor means that $J$ is
injective, while the last term being an exact functor means that
$J^\wedge$ is flat, since $\hhom(-,S^\wedge)$ is {\it a priori} a
faithful exact functor.  It follows that a complex $J^{\textstyle \cdot}
: 0 \to J^0 \to J^1 \to \cdots\ $ is an injective resolution of $N$ in
$\cal M$ if and only if $(J^{\textstyle \cdot})^\wedge$ is a flat
resolution of $N^\wedge$.  Substituting $k$ for $(-)$ in
Equation~(\ref{eqn:dash}) and applying Proposition~\ref{prop:calculate}
we get
\begin{equation} \label{eqn:ext.tor}
  \eext^i(k,N)\ \cong\ \ttor_i(k,N^\wedge)^\wedge
\end{equation}
from which the result follows at once.
\hfill
$\Box$
\begin{cor} \label{cor:lattice}
$\mu_{i,\bb}(T/M) = \beta_{i,-\bb}(M^T)$ for any monomial module $M
\subseteq T$.  \hfill $\Box$
\end{cor}

The next few results are preliminary to the theorems relating the Betti
numbers of $I$ to the Bass numbers of $I$ (Theorem~\ref{thm:gor}) and the
Bass numbers of $I^{[\aa]}$ (Theorem~\ref{thm:bass-betti}).
\begin{prop} \label{prop:Ibetti}
Let $I$ be an ideal.  Then
$$
\begin{array}{ll}
(i)
&	\beta_{i,\bb\,}(\,\widetilde{I}\,) = 0\,
	{\rm \ if\ } \bb \not\preceq \aa_I. \cr
(ii)
&	\beta_{i,\bb\,}(\,\widetilde{I}\,) = 0\, 
	{\rm \ if\ } \bb \not\succeq {\bf 1}. \cr
(iii)
&	\beta_{i,\bb\,}(\,\widetilde{I}\,) = \beta_{i,\bb\,}(\,I\,)\, 
	{\rm \ if\ } \,{\bf 1} \preceq \bb. \cr
\end{array}
$$
\end{prop}
{\it Proof:\ } Given any submodule $M \subseteq T$, define for each $\bb
\in \ZZ^n$ the following simplicial subcomplex of $\Delta$:
$$
  K_\bb(M) \ = \ \{F \in \Delta \mid x^{\bb - F} \in M\}\,.
$$
It is a result of \cite{Hoc} and \cite{Ros} (and extended to the case $M
\subseteq T$ by \cite{BS}) that
$$
  \beta_{i,\bb}(M) = \dim_k\HH_i(K_\bb(M);k), 
$$
the dimension of the $i^{\rm th}$ simplicial homology of $K_\bb(M)$ with
coefficients in $k$.  To prove (i) and (ii) it suffices to show that
$K_\bb(\widetilde{I})$ is a cone (and therefore acyclic) unless ${\bf 1}
\preceq \bb \preceq \aa_I$.  If $\aa_I = (a_1,\ldots,a_n)$ and $b_i \geq
a_i + 1$, then it follows from Proposition~\ref{prop:tilde}(ii) that
$K_\bb(\widetilde{I})$ is a cone with vertex $\{i\}$, proving (i).  That
$K_\bb(\widetilde{I})$ is a cone with vertex $\{i\}$ if $b_i \leq 0$
follows directly from the definition of \v Cech hull, proving (ii).
Finally, (iii) holds because $K_\bb(\widetilde{I}) = K_\bb(I)$ whenever
$\bb \succeq {\bf 1}$.
\hfill
$\Box$
\begin{lemma} \label{lemma:cutoff}
Let $M \subseteq T$.  Then $\beta_{i,\bb}(\,M\,) = \beta_{i,\bb}(\,M / M
\cap \widetilde{\mm^{\aa + {\bf 1}}}\,)$ if \,$\bb \preceq \aa$.
\end{lemma}
{\it Proof:\ } It follows from Example~\ref{ex:tilde}(ii) that $(M \cap
\widetilde{\mm^{\aa + {\bf 1}}})_\bb = 0$ if \,$\bb \preceq \aa$, so the
Taylor resolution of it (see \cite{Tay} for the original or \cite{BS},
Proposition~1.5 for a treatment including submodules of $T$) forces
$\beta_{i,\bb}(\,M \cap \widetilde{\mm^{\aa + {\bf 1}}}\,) = 0$ for all
$\bb \preceq \aa$.  Applying $\ttor$ to the exact sequence
$$
  0 \to M \cap \widetilde{\mm^{\aa + {\bf 1}}} \to M \to M / M \cap
  \widetilde{\mm^{\aa + {\bf 1}}}
$$
yields the result.
\hfill
$\Box$

\begin{lemma} \label{lemma:mu_n}
If $i < n$ then $\eext^i(k,S/I) \cong \eext^i(k,T/I)$, and in the
remaining case $i = n$ we have $\eext^n(k,S/I) = k[{\bf 1}]$.
\end{lemma}
{\it Proof:\ } One can first calculate 
$
\eext^i(k,S) = 
	\left\{\begin{array}{ll}
		k[{\bf 1}] & {\rm if}\ i = n \cr
		0	   & {\rm otherwise}
	\end{array}\right.
$
from the Koszul complex and $\eext^i(k,T) = 0$ for all $i$ because $T$ is
injective in the category $\cal M$.  Using the long exact sequence of
$\eext$ from $\,0 \to S \to T \to T/S \to 0\,$ one finds that
$\eext^i(k,S) \cong \eext^{i-1}(k,T/S)$.

{}From the above calculations and the long exact sequence of $\eext$
arising from
$$
	0 \to S/I \to T/I \to T/S \to 0
$$
the lemma will follow if we can show that the map
$$
	\eext^{n-1}(k,T/S) \to \eext^n(k,S/I)
$$
is an isomorphism.  But $S$ is a regular ring, so $\eext^n(k,S/I)$ is
nonzero {\it a priori} because of \cite{BH}, Proposition~3.1.14 and
\cite{BH}, Theorem~3.1.17, so it is enough to prove that the
$1$-dimensional vector space $\eext^{n-1}(k,T/S) \cong \eext^n(k,S) \cong
k[{\bf 1}]$ maps surjectively, i.e.\ that $\eext^n(k,T/I) = 0$.  Now
$\eext^n(k,T/I) \cong \eext^{n+1}(k,I)$ because of the exact sequence
$$
	0 \to I \to T \to T/I \to 0,
$$
and $\eext^{n+1}(k,I) = 0$ because of the same \cite{BH} reference as
above.
\hfill
$\Box$
\vskip 2mm

The next main result, Theorem~\ref{thm:gor}, is really a rephrasing of an
observation made in the proof of \cite{Hoc}, Theorem~5.2.  While it is
possible, by quoting the self-duality of the Koszul complex, to extend
the result to include all $S$-modules, the proof here demonstrates
effectively the interaction of Alexander duality with other kinds of
duality.  Aside from the intrinsic interest in its proof,
Theorem~\ref{thm:gor} will find an application in the proof of
Theorem~\ref{thm:bass-betti}.  Two preliminary results are needed, the
first of which will also be used in the proof of
Proposition~\ref{prop:restrict-localize}.
\begin{lemma} \label{lemma:Jtilde}
With $J = I + \mm^{\aa + {\bf 1}}$ we have $\widetilde{J}^{\;T} =\,
I^{[\aa]}[\aa]$.  The same is true if $I$ and $I^{[\aa]}$ are reversed.
\end{lemma}
{\it Proof:\ } The last statement is because of
Proposition~\ref{prop:Iaa=I}.  By Example~\ref{ex:tilde}(ii) and
Proposition~\ref{prop:T-dual},
\vskip 1mm
\noindent
\parbox{.95\linewidth}{\centering
$ \displaystyle
  I^{[\aa]}[\aa]\ =\ \widetilde{I}^{\;T} \cap S[\aa]\ =\ (\widetilde{I} +
  S[\aa]^{\:T\,})^T\ =\ \widetilde{J}^{\;T}.
$
}
{
\settodepth{\Mylen}{
\parbox{.95\linewidth}{\centering
$ \displaystyle
  I^{[\aa]}[\aa]\ =\ \widetilde{I}^{\;T} \cap S[\aa]\ =\ (\widetilde{I} +
  S[\aa]^{\:T\,})^T\ =\ \widetilde{J}^{\;T}.
$
}}
\addtolength{\Mylen}{-1.3mm}
\hfill \raisebox{-\Mylen}{$\Box$}
\vskip 1mm
}

The reader knowledgeable about linkage will recognize a hint of
Theorem~\ref{thm:linkage} in the next proposition.  Only the special case
$\bb = {\bf 0}$ is required in this section.  However, the more general
result is a major component in the proof of
Theorem~\ref{thm:coresolution}.
\begin{prop} \label{prop:ext}
Let $\aa \succeq \aa_I$, $J = I^{[\aa]} + \mm^{\aa + {\bf 1}}$, and $\bb
\in \NN^n$.  Then
$$
  \eext^n_S \Bigl( S[\bb]/ \,S[\bb]\!\cap\!\widetilde{J}, S \Bigr) \ = \
  \Bigl(I\,/\:I\!\cap\!\mm^{\aa + \bb + {\bf 1}}\Bigr)[\aa + {\bf 1}].
$$
In particular, taking $\bb = {\bf 0}$ yields $\eext^n(S/J,S) \cong
\Bigl(I/{\scriptstyle \,}I\!\cap\!\mm^{\aa + {\bf 1}} \Bigr) [\aa + {\bf
1}]$.
\end{prop}
{\it Proof:\ } The module $T/\widetilde{J}$ is the $k$-dual of the
finitely generated module $I[\aa]$ by Lemma~\ref{lemma:Jtilde}, and is
hence artinian by Matlis duality, cf.~\cite{GW}, Theorem~2.1.4.  Thus our
module $S[\bb]/ \,S[\bb]\!\cap\!\widetilde{J} \subseteq T/\widetilde{J}$
is also artinian, and (obviously) finitely generated, as well.  Since the
canonical module of $S$ is $S[-{\bf 1}]$ by \cite{GW}, Corollary~2.2.6,
local duality (in the form of \cite{GW}, Theorem~2.2.2) applied to the
zeroeth local cohomology module implies the first equality of the
following:
\newlength{\sumlen} \settowidth{\sumlen}{\mbox{$\displaystyle \sum_{.}$}}
$$
\begin{array}{@{}rclll@{}}
\eext^n_S \Bigl( S[\bb]/ \,S[\bb]\!\cap\!\widetilde{J}, S \Bigr)
& \!\! = \!\! &
	\Bigl(S[\bb]/\,S[\bb]\!\cap\!\widetilde{J}\Bigr)^\wedge[{\bf 1}] &
	\hbox{by local duality}
\phantom{\displaystyle \sum_{.}} \hbox{\hskip -\sumlen}
\cr
& \!\! = \!\! &
	\Bigl(\widetilde{J}^{\;T}/\,\widetilde{J}^{\;T}\!\!\cap\! 
	S[\bb]^T\Bigr)[{\bf 1}] & 
\phantom{\displaystyle \sum_{.}} \hbox{\hskip -\sumlen}
	\hbox{by Proposition~\ref{prop:T-dual}(vii)}
\cr
& \!\! = \!\! &
	\Bigl(I/\,I \!\cap\! S[\aa+\bb]^T\Bigr)[\aa + {\bf 1}] & 
\phantom{\displaystyle \sum_{.}} \hbox{\hskip -\sumlen}
	\hbox{by Lemma~\ref{lemma:Jtilde}, shifting by } [-\aa][\aa]
\cr
& \!\! = \!\! &
	\Bigl(I/\:I\!\cap\!\mm^{\aa + \bb + {\bf 1}}\Bigr)[\aa + {\bf 1}]
	&
	\hbox{by Example~\ref{ex:tilde}(ii)}. & \hskip -18pt \Box
\end{array}
$$

For artinian $J$, the list of Betti numbers for the canonical module
$\eext^n(S/J,S[-{\bf 1}])$ of $S/J$ is essentially the reverse of the
list of Betti numbers for $J$; see, for instance, \cite{BH},
Corollary~3.3.9.  On the other hand, there is the lattice-complementation
view of Alexander duality, which emerges in Corollary~\ref{cor:lattice}
as a relation between the Betti numbers of a monomial module and the Bass
numbers of its $T$-dual.  These two dualities can be combined to compare
the Betti numbers of $I$ to the Bass numbers of the same ideal $I$:
\begin{thm} \label{thm:gor}
For all $i \in \ZZ$ and $\bb \in \ZZ^n$,
$$
	\beta_{\,n-i\,,\,\bb\,}(\,S/I\,)\ =\ 
	\mu_{\,i\,,\,\bb - {\bf 1}\,}(\,S/I\,).
$$
\end{thm}
{\it Proof:\ } The case $i = n$ follows from the calculations of
Lemma~\ref{lemma:mu_n}, so assume from now on that $i \leq n-1$.  In
particular, we can calculate the Bass numbers from $T/I$ instead of $S/I$
by Lemma~\ref{lemma:mu_n}.  Let $\aa = \aa_I + {\bf 1}$.  All of the
Betti numbers are zero in degrees $\bb \not\preceq \aa$ by
Proposition~\ref{prop:Ibetti}(i) and~(iii).  As for the Bass numbers, we
can use the fact that, with $J := I^{[\aa]} + \mm^{\aa + {\bf 1}}$, we
have $I^{\,T} = \widetilde{J}[\aa]$ by Lemma~\ref{lemma:Jtilde}.  It
follows that $\mu_{i,\bb - {\bf 1}}(\,T/I\,) = \beta_{i,{\bf 1} -
\bb}(\,\widetilde{J}[\aa]\,) = \beta_{i,\aa + {\bf 1} -
\bb}(\,\widetilde{J}\,)$ by Corollary~\ref{cor:lattice}, and then
Proposition~\ref{prop:Ibetti}(ii) implies that these numbers are zero if
$\bb \not\preceq \aa$.

{}From now on, assume $\bb \preceq \aa$ and $0 \leq i \leq n-1$.  Let $J =
I^{[\aa]} + \mm^{\aa + {\bf 1}}$ and calculate
$$
\begin{array}{rcll}
\mu_{\,i\,,\,\bb - {\bf 1}\,}(\,S/I\,)
& = &
	\mu_{\,i\,,\,\bb - {\bf 1}\,}(\,T/I\,)
&	\hbox{by Lemma~\ref{lemma:mu_n} and } i \leq n-1
\cr
& = &
	\beta_{\,i\,,\,\aa + {\bf 1} - \bb\,}(\,\widetilde{J}\,)
&	\hbox{by Corollary~\ref{cor:lattice}, since } I^{\,T} =
		\widetilde{J}[\aa]
\cr
& = &
	\beta_{\,i\,,\,\aa + {\bf 1} - \bb\,}(\,J\,)
&	\hbox{by Proposition~\ref{prop:Ibetti}(iii) and }\bb \preceq \aa
\cr
& = &	\beta_{\,i+1\,,\,\aa + {\bf 1} - \bb\,}(\,S/J\,)
&	{\rm since}\ i \geq 0
\cr
& = &
	\beta_{\,n-i-1\,,\,\bb - {\bf 1} - \aa\,}
		\Bigl(\eext^n(S/J,S)\Bigr) 
&	\hbox{since $S$ is Gorenstein, $S/J$ artinian}
\cr
& = &
	\beta_{\,n-i-1\,,\,\bb\,}
		\Bigl((I/I\!\cap\!\mm^{\aa + {\bf 1}})\Bigr)
&	\hbox{by Proposition~\ref{prop:ext}}
\cr
& = &
	\beta_{\,n-i-1\,,\,\bb\,}(\,I\,)
&	\hbox{by Lemma~\ref{lemma:cutoff} and } \bb \preceq \aa
\cr
& = &
	\beta_{\,n-i\,,\,\bb\,}(\,S/I\,)
&	\hbox{since } i \leq n-1
\cr
\end{array}
$$
proving the theorem.
\hfill
$\Box$

\end{section}
\begin{section}{Localization and restriction}

\label{section:localize}

This section aims to reveal the equality (Theorem~\ref{thm:bass-betti})
between Betti numbers of $I$ and localized Bass numbers
(Definition~\ref{defn:bass}) of $I^{[\aa]}$.  This equality generalizes
Theorem~\ref{thm:[a]=a}.  As a consequence of the equality, an inequality
between Betti numbers of $I$ and $I^{[\aa]}$ is obtained in
Theorem~\ref{thm:inequalities}, generalizing to arbitrary monomial ideals
an inequality of \cite{BCP} which was proven for radical ideals.

The next proposition should be thought of as the nonlocalized precursor
to Theorem~\ref{thm:bass-betti}(i).  As usual, let $\aa \geq \aa_I$.
\begin{prop} \label{prop:precursor}
If $I$ is an ideal and \,${\bf 1} \preceq \bb \preceq \aa$ then
$\beta_{i,\bb}(\,I\,) = \mu_{i,\bb^{\scriptstyle \aa} - {\bf 1}}(\,S /
I^{[\aa]}\,)$.
\end{prop}
{\it Proof:\ } If, to start with, $\bb \preceq \aa$, then 
$$
\begin{array}{@{}rcll@{}}
\beta_{i\,,\,\bb\,}(\,M\,)
& = &
	\mu_{i\,,\,-\bb\,}\Bigl( (M / M \cap \widetilde{\mm^{\aa + {\bf
	1}}})^\wedge \Bigr) &
	\hbox{by Lemma~\ref{lemma:cutoff} and
		Proposition~\ref{prop:complementation} }
\phantom{\displaystyle \sum_{.}} \cr
& = &
	\mu_{i\,,\,-\bb\,}\Bigl( S[\aa] / S[\aa] \cap M^T \Bigr) &
	\hbox{by Proposition~\ref{prop:T-dual}(vii) and 
	Example~\ref{ex:tilde}(ii)}
\phantom{\displaystyle \sum_{.}} \cr
& = &
	\mu_{i\,,\,\aa-\bb\,}\Bigl( S / M^T[-\aa] \cap S
	\Bigr). & \cr
\end{array}
$$
Substituting $M = \widetilde{I}$ we get $\beta_{i\,,\,\bb}
(\,\widetilde{I}\,) = \mu_{i\,,\,\aa-\bb}(\,S/I^{[\aa]}\,) $ if $\bb
\preceq \aa$, and when the assumption ${\bf 1} \preceq \bb$ is added,
$\aa-\bb = \bb^\aa - {\bf 1}$ and  the result is a consequence of
Proposition~\ref{prop:Ibetti}(iii).
\hfill
$\Box$
\vskip 2mm

Theorem~\ref{thm:bass-betti} is the combination of the previous
proposition with localization and restriction of scalars.  The following
definitions will provide concise notation for these operations, which
will be needed also for the definition of Bass numbers at primes other
than $\mm$ (Definition~\ref{defn:bass}).  Recall that $\overline{F} =
\{1,\ldots,n\} \setminus F = {\bf 1} - F$.
\begin{defn} \label{defn:localization}
Let $\Delta$ be the $(n-1)$-simplex on the vertices $\{1,\ldots,n\}$ and
$F \in \Delta$.  Define
$$
\begin{array}{lrcll}
(i)
& 	N(-\overline{F}^{\,}) & := & S[x^{-\overline{F}\:}] \otimes_S N
&	\hbox{for arbitrary modules}\ N \cr
(ii)
&	S_{[F]} & := & k[x_i \mid i \in F]
&	{\it a}\ \ZZ^F\!\hbox{-graded}\ k\hbox{-subalgebra of}\ S \cr
(iii)
&	N_{[F]} & := & \bigoplus_{\bb \in \ZZ^F} N_\bb
&	{\it a}\ \ZZ^F\!\hbox{-graded}\ S_{[F]}\hbox{-module} \cr
(iv)
&	N_{(F)} & := & N(-\overline{F}^{\,})_{[F]} \cr
\end{array}
$$
\end{defn}
The operations on $N$ listed above are all exact and commute with sums.
They should be thought of as: (i) homogeneous localization at $\mm^{F}$,
(iii) taking the ``degree zero part'' of $N$ with respect to $F$, and
(iv) taking the ``degree zero part of the homogeneous localization at
$\mm^{F\,}$'' as in algebraic geometry.  In (ii) and (iii), the copy of
$\ZZ^F$ may be thought of as sitting inside $\ZZ^n$ in the obvious way:
as the space spanned by the basis vectors in the support of $F$.
Thinking of $\ZZ^F$ this way can cause notational problems, however.  For
instance, any $\ZZ^n$-graded $S$-module $N$ can be thought of as a
$\ZZ^F$-graded $S_{[F]}$-module which in degree $\bb \in \ZZ^F$ is
$$
	\bigoplus_{\cc \cdot F\,=\,\bb} N_\cc \ =
	\bigoplus_{\bb'\,\in\,\ZZ^{\overline F}} N_{\bb\,+\,\bb'}\,,
$$
where $\cc \cdot F$ denotes the restriction to $F$ as in
Section~\ref{section:definitions}.  Note that the right-hand side gives
this vector space the structure of a $\ZZ^{\overline F}$-graded
$S_{[\,{\overline F}\,]}$-module.  The convention will be the following:

\noindent
{\narrower\vbox{\noindent
If $N$ is a $\ZZ^F$-graded $S_{[F]}$-module and $\bb \in \ZZ^F$, the
graded piece of $N$ in degree $\bb$ will be denoted $N_{\bb \cdot F}$.
That way, if $N$ happens also to be a $\ZZ^n$-graded $S$-module, the
usual notation $N_\bb$ can retain its old meaning as the degree $\bb$
part in the $\ZZ^n$-grading.
}}

\noindent
Even if $\bb \not\in \ZZ^F$ it will sometimes be convenient to use
$N_{\bb \cdot F}$ to denote the $\bb \cdot F$ graded piece in the
$\ZZ^F$-grading; i.e.\ with $\cc = \bb \cdot F \in \ZZ^F$, we set $N_{\bb
\cdot F} := N_{\cc \cdot F}$.  The next Lemma follows from the
definitions and the convention above.  In each of (i)--(v), the objects
are $\ZZ^F$-graded $S_{[F\,]}$-modules, but in (i), the objects may also
be considered as $\ZZ^{\overline F}$-graded $S_{[\,{\overline
F}\,]}$-modules or even $\ZZ^{\overline F} \times \ZZ^F = \ZZ^n$-graded
$S_{[\,{\overline F}\,]} \otimes_k S_{[F\,]} = S$-modules.

\begin{lemma} \label{lemma:localization}
For any $F \in \Delta$,
$$
\begin{array}{lrcl}
(i)
&	M(-\overline{F}^{\,}) & = & T_{[\overline{F}]} \otimes_k M_{(F)}
\ =\ 	S(-\overline{F}^{\,}) \otimes_{S_{(F)}} M_{(F)}
\phantom{\Bigl(\Bigr)}
\cr
(ii)
&	M_{[F]} & = & M_{{\bf 0} \cdot \overline{F}}
\phantom{\Bigl(\Bigr)}
\cr
(iii)
&	M[\aa]_{[F]} & = & M_{\aa \cdot \overline{F}\,}[\aa \cdot F]
\phantom{\biggl(\biggr)}
\cr
(iv)
&	(\widetilde{I^{\,}})_{[F]} & = & \widetilde{I_{\,[F]}}
\phantom{\Bigl(\Bigr)}
\cr
(v)
& 	(M^{\,T\,})_{[F]} & = & (M_{[F]})^{\,T_{[F]}}
\phantom{\Bigl(\Bigr)}
\cr
\end{array}
$$
where the right-hand sides of (iv) and (v) are, respectively, the \v Cech
hull and $T$-dual over $S_{[F]}$. \hfill $\Box$
\end{lemma} 

For submodules $M \subseteq T$ the various gradings allow for convenient
characterizations of localiza\-tion as in
Definition~\ref{defn:localization}(iv).  They use the fact that for any
$\bb \in \ZZ^n$, $M_{\bb \cdot \overline{F}}\,$ is naturally a submodule
of $T_{[F]} = T_{(F)}$.

\begin{prop} \label{prop:unions}
Let $M$ be a monomial module.
$$
\begin{array}{lrcl}
(i)\phantom{i}
&	M_{(F)}
& = &
	{\displaystyle \bigcup_{\bb\,\in\,\ZZ^{\overline{F}}}}
	M_{\bb \cdot \overline{F}}\:.
\end{array}
$$
If $M$ can be generated in degrees $\cc$ satisfying $\cc \cdot
\overline{F} \preceq \aa \cdot \overline{F}$ then
$$
\begin{array}{lrcl}
(ii)
&	M_{(F)}
& = &
	M_{\aa \cdot \overline{F}}\:.
	\phantom{\scriptstyle \bb\,\in\,\ZZ^{\overline{F}}\!}
\end{array}
$$
\end{prop}
{\it Proof:\ } (i) Observe that $M \subseteq M(-\overline{F}^{\,})
\subseteq T$ because everything is torsion-free.  Thus, if $\bb \in
\ZZ^{\overline{F}}$, then multiplication by $x^{-\bb}$ induces an
inclusion $M_{\bb \cdot \overline{F}} \to M_{(F)}$.  For the other
inclusion, note that any monomial in $M_{(F)}$ can be written as $x^\bb
\cdot x^\cc$ for some $x^\cc \in M$ and $\bb = -(\cc \cdot \overline{F})
\in \ZZ^{\overline{F}}$.

\vskip 1mm \noindent
(ii) The collection $\{M_{\bb \cdot \overline{F}}\}_{\bb\,\in\,
\ZZ^{\overline F}}\,$ of $S_{[F]}$-submodules of $T_{[F]}$ is partially
ordered by inclusion because $M$ is a module.  The union in (i)
stabilizes after $\aa \cdot \overline{F}$ if $M$ is generated in degrees
$\cc$ satisfying $\cc \cdot \overline{F} \preceq \aa \cdot \overline{F}$.
\hfill
$\Box$
\begin{example}\rm
Figure~3 illustrates some parts of Definition~\ref{defn:localization} and
Lemma~\ref{lemma:localization} in a specific case.  For notation, $x$,
$y$, and $z$ are identified with $1$, $2$, and $3 \in \{1,\ldots,3\} =
\Delta$.  The face $F$ is $\{y,z\} = \{2,3\}$, while $\bb = (4,4,2)$.
The small colored dots represent generators or irreducible components in
the restricted ideals.  It is not true that $\bb \succeq \aa_I$, so
Proposition~\ref{prop:unions} does not apply; nevertheless, $I_{\bb \cdot
{\overline F}} = I_{(F)}$ for these $\bb$, $I$, and $F$.  Figure~3 can
also be used as a test case for Proposition~\ref{prop:restrict-localize}.

\begin{figure}
\centering \includegraphics{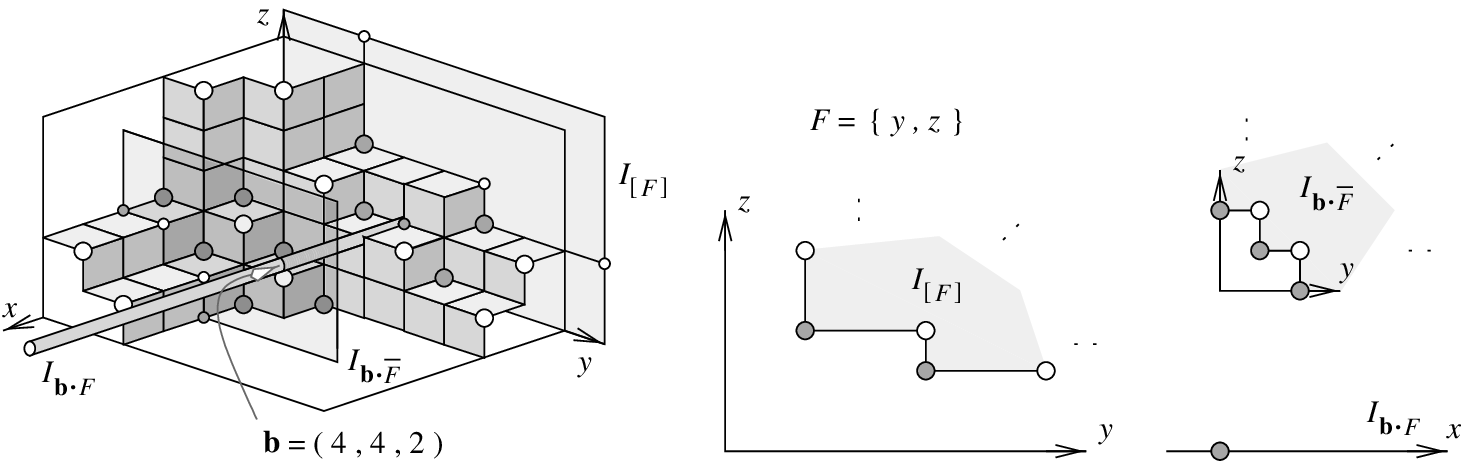}\\
Figure~3
\end{figure}
\end{example}
\begin{prop} \label{prop:restrict-localize}
$(I^{[\aa]})_{(F)} = (I_{\,[F]})^{[\aa \cdot F]}$ as ideals in $S_{(F)} =
S_{[F]}$.  In words, dualizing and then localizing is the same as
restricting and then dualizing.
\end{prop}
{\it Proof:\ } It is enough to show that $(I^{[\aa]})_{(F)}[\aa \cdot F]
= (I_{\,[F]})^{[\aa \cdot F]}[\aa \cdot F]$.  Now
$$
\begin{array}{rcll}
(I^{[\aa]})_{(F)}[\aa \cdot F]
& = &
	(I^{[\aa]})_{\aa \cdot \overline{F}}[\aa \cdot F]
	\phantom{\biggl(\biggr)}
&	\hbox{by Proposition~\ref{prop:unions}(ii) and
		Proposition~\ref{prop:Iaa=I}}
\cr
& = &
	(I^{[\aa]}[\aa])_{[F]}
&	\hbox{by Lemma~\ref{lemma:localization}(iii)}
\cr
& = &
	\Bigl(\Bigl((I + \mm^{\aa + {\bf 1}})\!\widetilde{\phantom{ii}}
	\Bigr)^T\Bigr)_{[F]}
&	\hbox{by Lemma~\ref{lemma:Jtilde},}
\cr
\end{array}
$$
and one can use the rules~\ref{lemma:localization}(v) and
then~\ref{lemma:localization}(iv) for interchanging the various
operations to get the last line to equal
$$
  \Bigl((I_{[F]} + \mm^{\aa \cdot F + F}{}_{[F]})\!\widetilde{\phantom{ii}}
  \Bigr)^{T_{[F]}},
$$
where $(-)^{T_{[F]}}$ is $T$-duality over $S_{[F]}$ as in
Lemma~\ref{lemma:localization}(v).  Another application of
Lemma~\ref{lemma:Jtilde} (over $S_{[F]}$ this time) gives the desired
result.
\hfill
$\Box$

\begin{prop} \label{prop:restriction}
Let $I \subseteq S$ and $\bb \in \ZZ^F$.  Then $\beta_{i,\bb}(I^{\,}) =
\beta_{i,\bb \cdot F}(I_{[F]})$.
\end{prop}
{\it Proof:\ } Let $\FF$ be the Taylor resolution of $I$ (see the proof
of Lemma~\ref{lemma:cutoff} for references).  Then $\FF_{[F]}$ is the
Taylor resolution of $I_{[F]}$.  Furthermore, $(k \otimes_S \FF)_{[F]} =
k \otimes_{S_{[F]}} \FF_{[F]}$ because if $\bb \in \NN^n$ then
$$
  \Bigl(k \otimes_S S[-\bb]\Bigr)_{[F]}
  \ = \ 
  k \otimes_{S_{[F]}} S[-\bb]_{[F]}
  \  = \ 
  \left\{\begin{array}{ll}
	k[-\bb]	& {\rm if}\ \bb \in \ZZ^F \cr
	0	& {\rm if}\ \bb \not\in \ZZ^F \cr
	\end{array}\right..
$$
Thus the desired Betti numbers are calculated from the same complex of
$k$-vector spaces.
\hfill
$\Box$

\begin{defn}[Bass numbers for arbitrary monomial primes] \label{defn:bass}
Given a module $N$ and a degree $\bb \in \ZZ^F$, the \emph{$i^{\rm \,th}$
Bass number of $N$ with respect to $F$ (or the prime ideal $\mm^{F}$) in
degree $\bb$} is defined as
$$
  \mu_{i,\bb}(F,N)\ :=\ \dim_k \Bigl(
  \eext_{S_{(F)}}^i(k,N_{(F)})_\bb \Bigr)\,.
$$
\end{defn}
\begin{remark} \rm
When $F = {\bf 1}$ this definition agrees with the Bass numbers of
Definition~\ref{defn:bass-betti}.
\end{remark}

Now comes the main result of this section.  It can be thought of as a
far-reaching generalization of Theorem~\ref{thm:[a]=a}, which is a
special case, pending the appropriate interpretation of Bass
numbers---see Proposition~\ref{prop:zeroeth_bass} and the second proof of
Theorem~\ref{thm:[a]=a} that follows it.  In part~(i) of the next
theorem, the case where $\bb$ has full support is
Proposition~\ref{prop:precursor}.
\begin{thm} \label{thm:bass-betti}
If $\:{\bf 0} \neq F \preceq \bb \preceq \aa \!\cdot\! F$ then for all $i
\in \ZZ$ we have
$$
\begin{array}{lrcl}
(i)
&	\beta_{i\,,\,\bb}(I) 
& = &
	\mu_{i\,,\,\bb^{\scriptstyle \aa} - F}(F,S/I^{[\aa]})
\cr
(ii)
&	\mu_{n-i-1\,,\,\bb - {\bf 1}}(S/I)
& = &
	\mu_{i\,,\,\bb^{\scriptstyle \aa} - F}(F,S/I^{[\aa]})
\cr
(iii)
&	\beta_{i\,,\,\bb}(I)
& = &
	\beta_{|F|-i-1\,,\,\bb^{\scriptstyle \aa}}
	(I^{[\aa]}{}_{(F)}).
\cr
\end{array}
$$
In any of these formulas, $I$ and $I^{[\aa]}$ can be switched, and the
same goes for $\bb$ and $\bb^\aa$.
\end{thm}
{\it Proof:\ } Statements (ii) and (iii) follow easily from (i), in view
of Theorem~\ref{thm:gor} and the fact that $\beta_{i,\bb}(I) =
\beta_{i+1,\bb}(S/I)$ when $\bb \neq {\bf 0}$.  To prove (i), note that
$\bb^\aa = (\bb \cdot F)^{\aa \cdot F}$, so
$$
\begin{array}{rcll}
\beta_{\,i\,,\,\bb\,}(\,I\,)
& = &
	\beta_{\,i\,,\,\bb \cdot F\,}(\,I_{\,[F]}\,)
	\phantom{\Bigl(\Bigr)}
&	\hbox{by Proposition~\ref{prop:restriction}}
\cr
& = &
	\mu_{\,i\,,\,\bb^{\scriptstyle \aa} - F\,}
		(\,S_{[F]} / I_{\,[F]}{}^{[\aa \cdot F]}\,) 
	\phantom{\Bigl(\Bigr)}
&	\hbox{by Proposition~\ref{prop:precursor}}
\cr
& = &
	\mu_{\,i\,,\,\bb^{\scriptstyle \aa} - F\,}
		(\,S_{(F)} / I^{[\aa]}{}_{\,(F)}\,) 
	\phantom{\Bigl(\Bigr)}
&	\hbox{by Proposition~\ref{prop:restrict-localize}}
\cr
& = &
	\mu_{\,i\,,\,\bb^{\scriptstyle \aa} - F\,}
		(\,F\,,\,S / I^{[\aa]}\,) 
	\phantom{\Bigl(\Bigr)}
&	\hbox{by definition}
\cr
\end{array}
$$
since $(-)_{(F)}$ is exact.  Note that the Bass number in the penultimate
line is with respect to the maximal ideal of $S_{(F)}$.  The last
statement of the theorem is true because $(\bb^\aa)^\aa = \bb$ and
$(I^{[\aa]})^{[\aa]} = I$, and because imposing the condition on $\bb$ is
equivalent to imposing the same condition on $\bb^\aa$.
\hfill
$\Box$

\begin{remark} \label{rk:links} \rm
Part (i) of the theorem can be thought of as the generalization to
arbitrary monomial ideals of the formulas in \cite{ER}, Proposition~1 and
\cite{BCP}, Theorem~2.4, using \cite{Hoc}, Theorem~5.2 and the fact that
links come from localization (\cite{Hoc}, Proposition~5.6).
\end{remark}

As a consequence of the theorem, the list of Betti numbers of $I^\aa$
will be independent of $\aa$, though the $\ZZ^n$-degrees in which they
occur will vary with $\aa$.  Indeed, the list of Betti numbers of $I^\aa$
is just the list of (localized) Bass numbers of $I$ by part~(i) of the
theorem.  Thus the collection of ideals that are dual to $I$ are very
closely related homologically.  This will be highlighted again in
Section~\ref{section:cellular} in terms of various geometrically defined
resolutions.

Before Remark~\ref{rk:links}, the symbol $I^\aa$ had not appeared in this
section (or the last) without brackets on the $\aa$; that is, none of the
results have been logically dependent on Definition~\ref{defn:alexdual}
or Theorem~\ref{thm:[a]=a}.  Therefore, Theorem~\ref{thm:bass-betti} can
be used to give a second proof of Theorem~\ref{thm:[a]=a}.  In fact, this
``second proof'' was discovered before the more elementary proof in
Section~\ref{section:alternate}.  The next proposition is what allows the
irreducible decomposition to be read off of the zeroeth Bass numbers just
as the minimal generators are read off the zeroeth Betti numbers.
\begin{prop} \label{prop:zeroeth_bass}
Given an ideal $I \subseteq S$ the following are equivalent for $\bb \in
\ZZ^F$: 
$$
\begin{array}{ll}
(i)
&	\mm^{\bb}{\rm \ is\ an\ irredundant\ irreducible\
	component\ of\ } I. \cr 
(ii)
&	\mu_{0,\bb - F}(F,S/I) = 1. \cr
(iii)
&	\mu_{0,\bb - F}(F,S/I) \neq 0. \cr
\end{array}
$$
\end{prop}
{\it Proof:\ } Let $I = \bigcap_j Q_j$ be the (unique) irredundant
decomposition of $I$ into irreducible ideals $Q_j$.  Then we have an
injection $\ 0 \to S/I \to \bigoplus_j S/Q_j\ $ which, by the proofs of
\cite{Vas}, Propositions~3.16 and~3.17, induces an {\it isomorphism}
\begin{equation} \label{eqn:hom}
  \hhom_S (S/\mm^F\!,S/I)(-\overline{F}^{\,}) \to 
  \bigoplus_j \hhom_S (S/\mm^F\!,S/Q_j)(-\overline{F}^{\,});
\end{equation}
this is because the functor $\Delta_\pp(\cdot)_\pp$ in the \cite{Vas}
reference is easily seen to be Hom$_{\,}(R/\pp,\cdot)_\pp$ (so we can
take $\pp = \mm^F$).  Using Lemma~\ref{lemma:localization}(i) we can move
the localization into and out of the $\hhom$: for any finitely generated
$S$-modules $L$ and $N$,
$$
\begin{array}{r@{\ \ \protect\cong\ \ }l}
\hhom_S \Bigl( L,N \Bigr)(-{\overline F})
&	\hhom_{S(-{\overline F})}
		\Bigl( L(-{\overline F}),N(-{\overline F}) \Bigr)
\cr
&	S(-{\overline F}) \otimes_{S_{(F)}} \hhom_{S_{(F)}}
		\Bigl( L_{(F)},N_{(F)} \Bigr)
	\phantom{\biggl(\biggr)}
\cr
&	T_{[\,{\overline F}\,]} \otimes_k \hhom_{S_{(F)}}
		\Bigl( L_{(F)},N_{(F)} \Bigr)\,.
\end{array}
$$
Treating these as $\ZZ^{\overline F}$-graded $S_{[\,{\overline
F}\,]}$-modules and taking the degree ${\bf 0} \!\cdot\! {\overline F}$
part in the last line yields $\hhom_{S_{(F)}} \Bigl( L_{(F)},N_{(F)}
\Bigr)$.  Applying this procedure to Equation~(\ref{eqn:hom}) reveals an
isomorphism
$$
  \hhom_{S_{(F)}} \Bigl(k,(S/I)_{(F)}\Bigr)\ \cong\ 
  \bigoplus_j \hhom_{S_{(F)}} \Bigl(k,(S/Q_j)_{(F)}\Bigr).
$$
Since we can calculate
$$
  \mu_{\,0,\bb - F\,}(F,S/Q_j)\ =\ \left\{\begin{array}{ll}
		1 & {\rm if}\ Q_j = \mm^\bb \cr
		0 & {\rm otherwise} \end{array}
	\right.
$$
the proposition follows from the definition of Bass numbers.
\hfill
$\Box$
\vskip 2mm

\noindent
{\it Second proof of Theorem~\ref{thm:[a]=a}:\ } Every generator $x^\bb$
of $I$ corresponds to a nontrivial zeroeth Betti number of $I$ which
satisfies the condition $F \preceq \bb \preceq \aa \cdot F$ for $F =
\sqrt\bb$ because $I \subseteq S$ and $\aa \succeq \aa_I$.  After
applying Theorem~\ref{thm:bass-betti}(i) and the previous proposition, we
can conclude that each generator of $I$ does indicate the presence of an
appropriate irreducible component of $I^{[\aa]}$.  To show that each
nontrivial zeroeth Bass number of $I^{[\aa]}$ comes from some Betti
number of $I$, we demonstrate that if $\bb \in \ZZ^F$ and $\mu_{0,\bb -
F}(F,S/I) \neq 0$ then $F \preceq \bb \preceq \aa \cdot F$.  Localizing
at $\mm^F$, we may assume that $F = {\bf 1}$.  Clearly $\bb \succeq {\bf
1}$ since $\mm^\bb$ is $\mm$-primary, so the desired result falls out of
Theorem~\ref{thm:gor} and Proposition~\ref{prop:Ibetti}.
\hfill
$\Box$
\vskip 2mm

Next on the agenda is the generalization to arbitrary monomial ideals of
an inequality of \cite{BCP} for squarefree ideals.  The topological
argument involving links employed there is preempted here by a simple
algebraic observation involving localization (which gives links in the
squarefree case, see \cite{Hoc}, Proposition~5.6).
\begin{thm} \label{thm:inequalities}
If $\:\aa \succeq \aa_I\:$ and $\:F \preceq \bb \preceq \aa \!\cdot\!
F\:$ then
$$
	\beta_{i\,,\,\bb}(I)\ \ \leq\ \sum_{\cc \cdot F \,=\,
	\bb^{\scriptstyle \aa}} \beta_{|F| - i - 1\,,\,\cc\,}(I^\aa).
$$
\end{thm}
{\it Proof:\ } Let $\FF$ be a minimal free resolution of $I^\aa$.
Localizing at $\mm^F$ we obtain a free resolution $\FF_{(F)}$ of
$I^\aa{}_{(F)}$ over $S_{(F)}$.  The generators of $\FF_{(F)}$ as a free
$S_{(F)}$-module are in bijective correspondence with the generators of
$\FF$ itself.  Specifically, for any $\bb' \in \ZZ^F$ we find that
$S[\cc]_{(F)} = S_{(F)}[\bb' \cdot F]$ if and only if $\cc \cdot F =
\bb'$.  Thus the number of summands of $\FF_{(F)}$ in homological degree
$|F| - i - 1$ and $\ZZ^F$-degree $\bb^\aa$ is equal to
$$
  \sum_{\cc \cdot F \,=\, \bb^{\scriptstyle \aa}} \beta_{|F| - i -
  1\,,\,\cc\,}(I^\aa)
$$
since $\FF$ is minimal.  On the other hand, the number of such summands
must clearly be $\geq \beta_{|F| - i - 1,\bb^{\scriptstyle \aa}}
(I^\aa{}_{(F)})$, with equality if and only if $\FF_{(F)}$ is minimal.
Since this last number is equal to $\beta_{i,\bb}(I)$ by
Theorem~\ref{thm:bass-betti}, we are done.
\hfill
$\Box$

\begin{cor}[Bayer-Charalambous-Popescu]
If $I$ is squarefree then 
$$
  \beta_{i\,,\,\bb}(I)\ \ \leq\ \sum_{\bb \preceq \cc \preceq {\bf 1}}
  \beta_{|\bb| - i - 1\,,\,\cc}(I^\vee)
$$
for $0 \leq i \leq n-1$ and ${\bf 0} \preceq \bb \preceq {\bf 1}$.
\end{cor}
{\it Proof:\ } This is a special case of the theorem once it is noted
that (i) $\beta_{|\bb| - i - 1, \cc}(I^\vee) = 0$ unless ${\bf 0} \preceq
\cc \preceq {\bf 1}$, and (ii) ${\bf 0} \preceq \cc$ and $\cc \cdot
\sqrt{\bb} = \bb$ imply $\cc \succeq \bb$.
\hfill
$\Box$

\end{section}
\begin{section}{Duality for cellular complexes: the cohull resolution}

\label{section:cellular}

\noindent
This section explores the effect of Alexander duality on various
geometrically defined free resolutions, in the spirit of \cite{BPS},
\cite{BS}, and \cite{Stu}.  First, the concept of a geometrically defined
resolution is broadened past {\it cellular resolutions} to include {\it
relative cocellular resolutions} (Definition~\ref{defn:relcocell}).  The
key result (Theorem~\ref{thm:relcocell}) is presented, though the
majority of its proof occupies Section~\ref{section:limits}.  As an
application, it is shown how irreducible decompositions can be specified
by cellular resolutions (Theorem~\ref{thm:irr}).  The culmination of
these ideas is a new canonical geometric resolution for monomial ideals
(Definition~\ref{defn:cohull}).  It is called the {\it cohull
resolution}, and is defined by applying Alexander duality to the hull
resolution of \cite{BS}.  As a special case, the co-Scarf resolution of a
cogeneric monomial ideal of \cite{Stu} is seen to be the cohull
resolution (Theorem~\ref{thm:coscarf}), and is thus Alexander dual to the
Scarf resolution of a generic monomial ideal \cite{BPS}.  A number of
examples are presented, including permutahedron and tree ideals.

\noindent
{\it Conventions regarding grading and chain complexes:}

{\narrower\noindent
A chain complex of $S$-modules 
$$
  \FF:\ \ \cdots \to N_{i+1} \to N_i \to N_{i-1} \to \cdots,\ \ \ N_i\
  {\rm in\ homological\ degree}\ i,
$$
is viewed as a (homologically) $\ZZ$-graded $S$-module $\bigoplus N_i$
with a differential $\partial$ of degree $-1$.  If ``$\,[\aa]$'' is
attached to $\FF$ then each summand is to be shifted in its
$\ZZ^n$-grading to the left by $\aa$, while ``$\,(j)$'' indicates that
the homological degrees are to be shifted down by $j$, yielding the
notation
$$
\FF[\aa](j):\ \ \cdots \to N_{i+1}[\aa] \to N_i[\aa] \to N_{i-1}[\aa]
\to \cdots,\ \ \ N_i\ {\rm in\ homological\ degree}\ i-j.
$$
Here, $N[\aa]_\bb = N_{\aa+\bb}$ for any $S$-module $N$ by definition.
Taking the $S$-dual $\FF^\Ast := \hhom(\FF,S)$ changes $\partial$ to its
transpose $\delta$, and makes homological degrees into cohomological
degrees, which are the negatives of homological degrees:
$$
\begin{array}{rl}
\FF^\Ast:\ \ \cdots \from N_{i+1}^\Ast \from N_i^\Ast \from
	N_{i-1}^\Ast \from \cdots, &
N_i^\Ast\ \hbox{in homological degree}\ -i
\cr
& \phantom{N_i^\Ast} =\ \hbox{cohomological degree}\ i.
\end{array}
$$
}

Labelled cell complexes provide compact vessels for recording the
monomial entries in certain $\ZZ^n$-graded free resolutions of an ideal.
\cite{BS} introduces this notion in the context of monomial modules, but
attention is restricted to boundary operators of the cell complex.  The
definitions below extend the concept to include coboundary operators, as
well.  For the reader's convenience, the definition of a labelled regular
cell complex and the cellular free complex it determines is recalled
briefly below, although the reader is urged to consult \cite{BS},
Section~1 for more details.

Let $\Lambda \subseteq \ZZ^n$ be a set of vectors, and let $X$ be a {\it
regular cell complex} whose vertices are indexed by the elements of
$\Lambda$.  For $\cc,\cc' \in \ZZ^n$, define the {\it join $\cc \vee
\cc'$} to be the componentwise maximum, i.e.\ $\cc \vee \cc'$ is the
smallest vector which is greater than or equal to both $\cc$ and $\cc'$
in each coordinate: $(\cc \vee \cc')_i = \max(c_i,c_i')$.  Given a face
$F \in X$, define the {\it label\, $\aa_F$ of $F$} to be the join
$\bigvee_{v \in F} \aa_v$ of the labels on the vertices in $F$, where
$\aa_v \in \Lambda$ is the element corresponding to $v$.  To avoid
confusion, the symbol $|X|$ will be used to denote the unlabelled
underlying cell complex of the labelled cell complex $X$.

We assume that $|X|$ comes equipped with an {\it incidence function}
$\varepsilon(F,F') \in \{1,0,-1\}$ defined on pairs of faces, which is
used to define the boundary map in the oriented chain complex of $|X|$
(with coefficients in $k$).  For each $F \in X$, let $SF$ be the free
$S$-module with one generator $F$ in degree $\aa_F$.  The {\it cellular
complex $\FF_X$} is the homologically and $\ZZ^n$-graded chain complex of
$S$-modules 

$$
  \FF_X \ = \bigoplus_{F \in X,\: F \neq \emptyset}SF
  \qquad {\rm with\ differential} \qquad
  \partial F \ = \sum_{G \in X,\,G \neq \emptyset}
		\varepsilon(F,G){m_F \over m_G} \cdot G,
$$
where $m_F := x^{\aa_F}$.  The homological degree of the basis vector $F
\in SF$ is the dimension of $F \in |X|$.  If $\FF_X$ is acyclic, it will
be said that $X$ {\it supports a free resolution} of the module
$\<x^{\aa_v} \mid v \in X$ is a vertex $\>$.

\noindent
\begin{remark} \rm
In \cite{BS} it is assumed that the elements of $\Lambda$ are pairwise
incomparable (as elements in the poset $\ZZ^n$), but $\Lambda$ is not
assumed to be finite.  Here, however, $\Lambda$ will always be finite,
but pairwise incomparability is not assumed.  It is easily verified that
all of the results in \cite{BS}, Section~1 remain true under these
hypotheses.
\end{remark}

\begin{defn}[Relative Cellular Complexes]
A \emph{relative cellular complex} $\FF_{\,(X,Y)}$ is the quotient
of a cellular complex $\FF_X$ supported on a labelled regular cell
complex $X$ by a subcomplex $\FF_Y$ for some regular cell
subcomplex $Y \subseteq X$, with the labelling on $Y$ induced by the
labelling on $X$.
\end{defn}

\begin{defn}[Relative Cocellular Complexes] \label{defn:relcocell}
A \emph{relative cocellular complex} $\FF^{\,(\hskip -1.2pt X,Y \hskip
-.4pt )}$ is obtained by taking $(\FF_{\,(X,Y)})^\Ast$ for a pair $(X,Y)$
of labelled relative regular cell complexes.  If $Y$ is empty,
$\FF^{\,(X,Y)}$ may be denoted $\FF^{\,X}\!$ and called a
\emph{cocellular complex} supported on $X$.
\end{defn}

\begin{remark} \rm
The relative cocellular complex $\FF^{\,(X,Y)}$ can be viewed as the
homogenization of the relative cochain complex of the pair $(X,Y)$, as
long as the label on a dual face $F^\Ast$ is the negative $-\aa_F$ of the
label on the face $F$.  The coboundary can then be written as $\delta
G^\Ast = \sum_{(F \in X,\,F \neq \emptyset)} \varepsilon(F,G){ m_F \over
m_G} \cdot F^\Ast$.
\end{remark}

\begin{defn}
Given a labelled regular cell complex $X$ and a vector $\bb \in \ZZ^n$,
define the following two labelled subcomplexes of $X$:
$$
\begin{array}{@{}l@{\ }l@{}}
\!(i)	&
	X_B(\bb) := \{F \in X \mid \aa_F \preceq \bb \}, \hbox{\rm\ the
	positively bounded subcomplex with respect to $\bb$}
\cr
\!(ii)\!&
	X_U(\bb) := \{F \in X \mid \aa_F \not\succeq \bb \}, 
	\hbox{\rm\ the negatively unbounded subcomplex with
	respect to $\bb$}
\end{array}
$$
Finally, let $X_U := X_U({\bf 1})$ be simply {\rm the negatively
unbounded subcomplex of $X$}.
\end{defn}
\begin{example}\rm 
\begin{figure}[htb] \label{fig:bounded}
\centering
\includegraphics{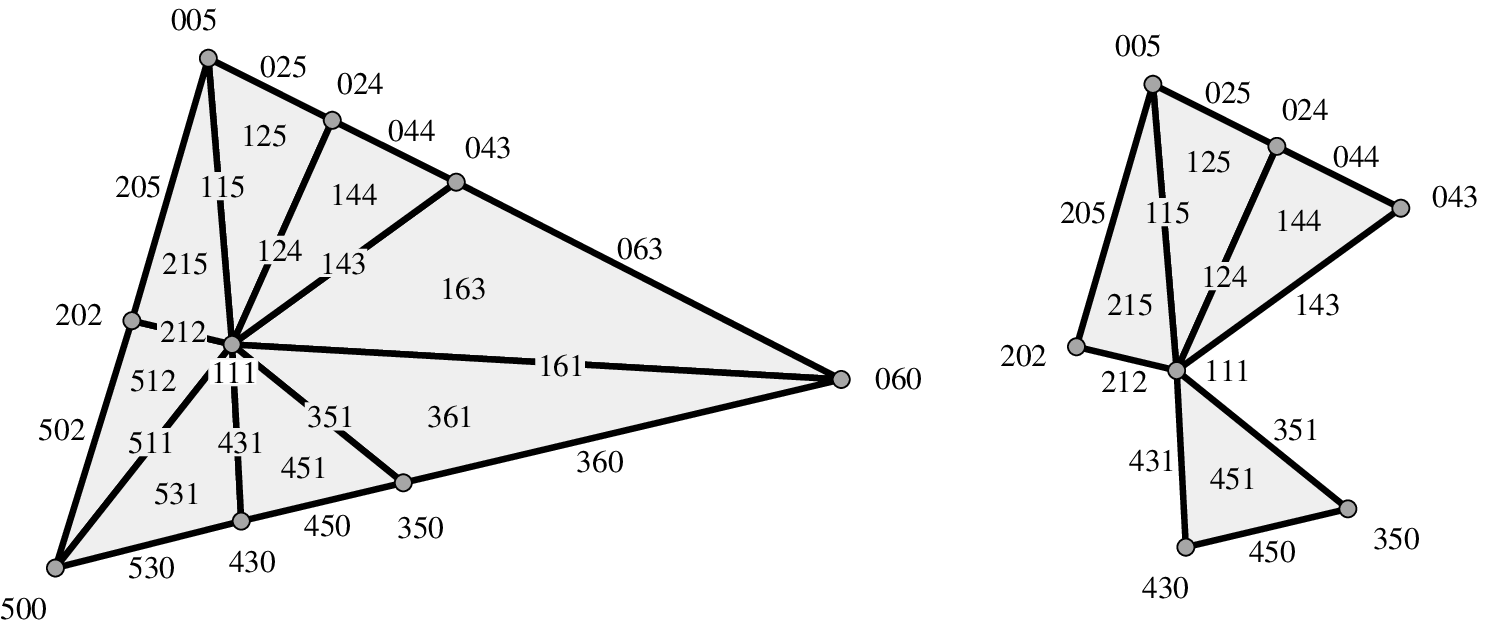}\\
\hskip 1in $X$ \hskip 2.8in $X_B(4,5,5)$\\
Figure~4
\end{figure}
Let $I$ be as in Example~\ref{ex:alexdual}.  The labelled complex $X$ in
Figure~4 is the {\it Scarf complex} \cite{BPS} of $I + \mm^{(5,6,6)}$
(see also Example~\ref{ex:scarf}, below).  Hence $\FF_X$ is a minimal
free resolution by \cite{BPS}, Theorem~3.2.  In this case, $(5,6,6) =
\aa_I + {\bf 1}$, but $z^5$ is already in $I$.  The label ``215'' in the
diagrams is short for (2,1,5).  The subcomplex $X_B(4,5,5)$, which is the
Scarf complex of the ideal $I$ itself, is also depicted in Figure~4 (see
Proposition~\ref{prop:cellres}, below).  The subcomplex $X_U$ is depicted
in Figure~5 along with a representation of the labelled relative cellular
complex $(X,X_U)$ and the relative cocellular complex $\FF^{(X,X_U)}$ of
free $S$-modules determined by it.
\begin{figure}[!htb] \label{fig:unbounded}
\centering
\includegraphics{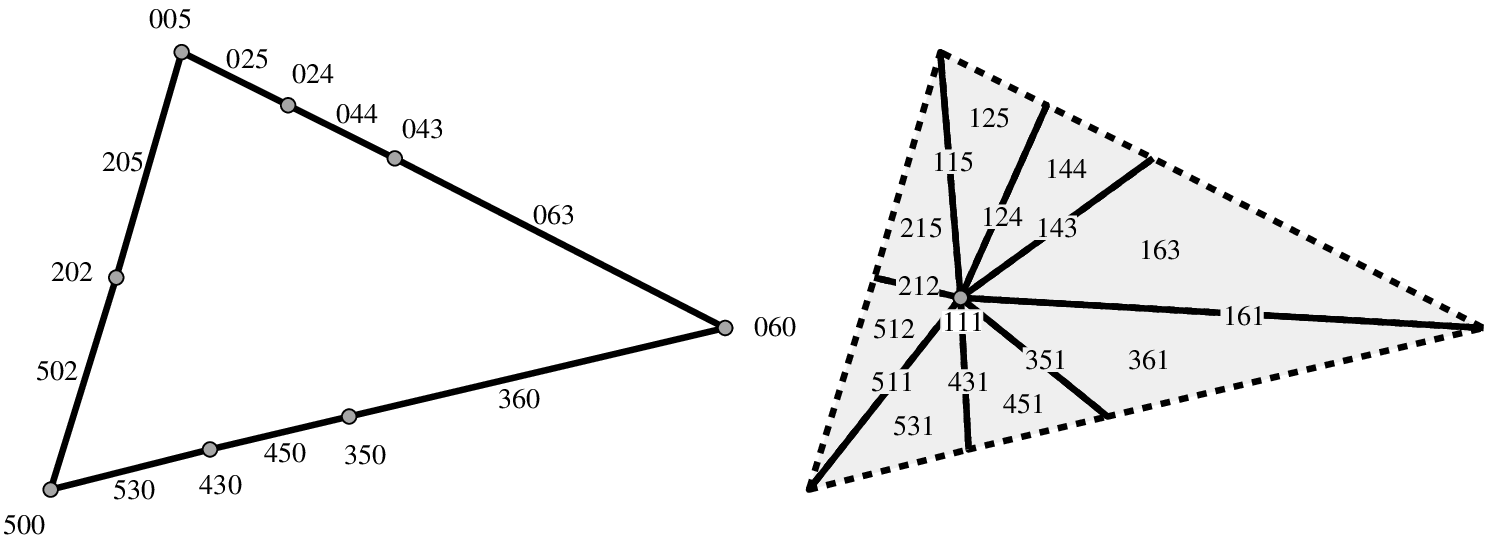}\\
\hfill $X_U = X_U({\bf 1})$ \hfill \hfill $(X,X_U)$ \hfill \mbox{}\\
$$
\begin{array}{@{}ccccccccc@{}} \scriptstyle
    &&&
	\left(\begin{array}{cccccccc}
		 x^2 & -y &   0  &   0  &   0  &  0  &  0  &  0  \cr
		  0  &  x & -y^2 &   0  &   0  &  0  &  0  &  0  \cr
		  0  &  0 &   x  & -y^2 &   0  &  0  &  0  &  0  \cr
		  0  &  0 &   0  &   z  & -x^3 &  0  &  0  &  0  \cr
		  0  &  0 &   0  &   0  &  z^3 & -x  &  0  &  0  \cr
		  0  &  0 &   0  &   0  &   0  &  y  & -z  &  0  \cr
		  0  &  0 &   0  &   0  &   0  &  0  & y^2 & -z  \cr
		-z^2 &  0 &   0  &   0  &   0  &  0  &  0  & y^2 \cr
	\end{array}\right) 
    &&
	\left(\begin{array}{c}
		y^5	\cr
		x^2y^4	\cr
		x^3y^2	\cr
		x^4	\cr
		xz	\cr
		z^4	\cr
		yz^3	\cr
		y^3z^2	\cr
	\end{array}\right)
    & 
    \cr
0 & \leftarrow & S^8 & \filleftmap & S^8 & \filleftmap & S & \leftarrow &
	0 \cr
\end{array}
$$
Figure~5
\end{figure}
For this, the edges have been oriented towards the center and the faces
counterclockwise.  The left copy of $S^8$ represents the 2-cells in
clockwise order starting from 361, the right copy of $S^8$ represents the
edges clockwise starting from 161, and the copy of $S$ represents the
lone vertex.  The other vertices and edges are not considered since they
lie in the subcomplex $X_U$.  It is not a coincidence that the negatively
unbounded subcomplex of $X$ is the topological boundary of $X$---this
will always happen for the Scarf complex of a generic artinian monomial
ideal, cf.\ Theorem~\ref{thm:boundary}.
\hfill
$\Box$
\end{example}

Recall that $\aa_I$ is the exponent on the least common multiple of the
minimal generators for $I$.  Suppose that we have a cellular resolution
$\FF_X$ of the ideal $I + \mm^{\aa + {\bf 1}}$ with $\aa \succeq \aa_I$.
\begin{prop} \label{prop:cellres}
$\FF_{X_B(\bb)}$ is a cellular resolution of $I$ for any $\bb$
such that $\aa_I \preceq \bb \preceq \aa$.
\end{prop}
{\it Proof:\ } With the conditions on $\bb$, the ideal $I$ is generated
by all monomials in $I + \mm^{\aa + {\bf 1}}$ whose exponent is $\preceq
\bb$, so the result is a direct consequence of \cite{BS}, Corollary~1.3.
\hfill $\Box$
\vskip 2mm

Duality for cellular resolutions says that if the cellular resolution
$\FF_X$ of $I +\mm^{\aa + {\bf 1}}$ has minimal length, a resolution
for the Alexander dual $I^\aa$ with respect to $\aa$ can also be
recovered from $X$:
\begin{thm} \label{thm:relcocell}
If the cellular resolution $\FF_X$ of $I + \mm^{\aa + {\bf 1}}$ has
length $n-1$ then $\FF^{\,(X,X_U)}[-\aa -{\bf 1}](1-n)$ is a relative
cocellular resolution of $I^\aa$.  Furthermore, this dual resolution is
minimal if $\FF_X$ is.
\end{thm}
{\it Proof:\ } The first statement will be a direct consequence of
Theorem~\ref{thm:coresolution}, below; the necessary assumption here that
$\FF_X$ has length $n-1$ is what makes $\FF^{\,(X,X_U)}[-\aa -{\bf
1}](1-n)$ a \emph{resolution} instead of just a free complex---that is,
there are no terms in negative homological degrees.  The construction of
$\FF^{\,(X,X_U)}$ from $\FF_X$ preserves minimality because the matrices
defining the differential of the former are submatrices of the transposes
of those defining the latter, and we need only check that these entries
are in $\mm$.
\hfill
$\Box$

\begin{remark} \label{rk:relcocell} \rm
(i) The hypothesis of the theorem requires that $X$ have dimension
$(n-1)$, so that $\FF_X$ has minimal {\it length}, but it does not
require that $\FF_X$ actually be a minimal resolution.

\vskip 1mm \noindent
(ii) It can be shown that $X_U$ may be replaced in the theorem by
$X_U(\bb+{\bf 1})$ for any $\bb$ satisfying ${\bf 0} \preceq \bb \preceq
\aa_I - \aa_{I^\vee}$.  Here, again, is the mysterious invariant from the
remark after Corollary~\ref{cor:Iaa=I}.  In most cases of interest,
though, $X_U = X_U(\bb)$ for all such $\bb$.
\end{remark}

\begin{example} \label{ex:coscarf} \rm
The free complex in Figure~5 is the minimal free resolution of the ideal
$I^\vee$ from Example~\ref{ex:alexdual}.  The reader may check, for
instance, that the product of the large matrix in Figure~5 with the list
of generators for $I^\vee$ (which may be treated as a matrix with one
row) is zero.  Note that the homological and $\ZZ^n$-graded shifts
promised by Theorem~\ref{thm:relcocell} aren't visible from the matrices.
\hfill $\Box$
\end{example}

Theorem~\ref{thm:relcocell} affords a generalization of \cite{BPS},
Theorem~3.7 on reading irreducible decompositions off of cellular
resolutions.  We will need the following.
\begin{lemma} \label{lemma:pure}
If the labelled cell complex $X$ supports a minimal free resolution of an
artinian ideal $J \subseteq S$ then $X$ is pure of dimension $n-1$.
\end{lemma}
{\it Proof:\ } Any facet has dimension $> 0$, so suppose that $F$ is a
facet of dimension $d > 0$.  Denote by $F^\Ast$ the basis element of the
cocellular complex $\FF^X$.  The modules $\eext\hdot(J,S)$ can be
calculated as the cohomology of $\FF^X$ by definition, and the coboundary
$\delta(F^\Ast)$ is zero because $F$ is a facet.  Moreover, the image of
$\delta$ is contained in $\mm\FF^X$ by minimality of $\FF_X$, whence
$F^\Ast$ is not itself a coboundary.  Thus $F^\Ast$ represents a nonzero
element of $\eext^d(J,S) \cong \eext^{d+1}(S/J,S)$.  It follows that $d =
n-1$ because $S/J$ has only one nonzero such $\eext$ module \cite{BH},
Proposition~3.3.3(b)(i).
\hfill
$\Box$
\vskip 2mm 

For the statement of the next theorem, the following notation is
convenient.  Suppose $\aa \succeq \aa_I$ and define, for any ${\bf 1}
\preceq \bb \preceq \aa + {\bf 1}$, the {\it bounded part} $\bb_B := (\aa
+ {\bf 1} - \bb)^\aa$ of $\bb$ with respect to $\aa$ to be the vector
whose $i^{\,\rm th}$ coordinate is $b_i$ if $b_i \leq a_i$ and zero if
$b_i = a_i + 1$.
\begin{thm} \label{thm:irr}
Let $\FF_X$ be a minimal cellular resolution of $I + \mm^{\aa + {\bf
1}}$.  Then the facets of $X$ are in bijection with the irredundant
irreducible components of $I$, and the intersection $\bigcap_F
\mm^{(\aa_F)_B}$ over all facets $F \in X$ is an irredundant irreducible
decomposition of $I$.
\end{thm}
{\it Proof:\ } It follows from Lemma~\ref{lemma:pure} that under these
conditions $X$ must be pure of dimension $n-1$.  Using this, it suffices
to show that the label on any facet is $\succeq {\bf 1}$, for then each
facet corresponds to a minimal generator of $I^\aa$ by
Theorem~\ref{thm:relcocell} and we are done by
Proposition~\ref{prop:min}.  Suppose, then, that $\aa_F$ is $\leq 0$ in
some coordinate for some facet $F$; say $(\aa_F)_n = 0$.  For $t >\!\!>
0$ consider $Y := X_B(t,t,\ldots,t,0)$, which gives a resolution of $J :=
(I + \mm^{\aa + {\bf 1}}) \cap k[x_1,\ldots,x_{n-1}]$ by \cite{BS},
Corollary~1.3.  The resolution $\FF_Y$ is minimal because $\FF_X$ is, and
$Y$ has dimension $n-1$ because $F \in Y$.  On the other hand, $J$ is an
artinian ideal of $k[x_1,\ldots,x_{n-1}]$, which contradicts
Lemma~\ref{lemma:pure} (with $n$ replaced by $n-1$).
\hfill
$\Box$
\vskip 2mm

The major consequence of Theorem~\ref{thm:relcocell} is the construction
of the cohull resolution (Definition~\ref{defn:cohull}) from the hull
resolution \cite{BS}, Section~2.  Therefore, we recall here the
definition of the hull complex. Let $t > (n+1)!$ and define $t^\bb :=
(t^{b_1},\ldots,t^{b_n})$.  The convex hull of the points $\{t^\bb \mid
x^\bb \in I\}$ is a polyhedron $P_t$ whose face poset is independent of
$t$.  It is shown in \cite{BS} that the vertices of $P_t$ are given by
those $t^\bb$ such that $x^\bb$ is a minimal generator of $I$.  The hull
complex hull$(I)$ is defined to be the bounded faces of $P_t$, but it may
also be described as those faces of $P_t$ admitting a strictly positive
inner normal.  The hull complex is labelled via the labels on its
vertices.
\begin{thm}[Bayer-Sturmfels] \label{thm:hull}
The free complex $\FF_{{\rm hull}(I)}$ is a cellular resolution of $I$.
\end{thm}
\begin{example} \label{ex:scarf} \rm
Let $\Lambda$ be the set of exponents on the minimal generators of a
generic monomial ideal $I$, and let $X$ be the labelled simplex with
vertices in $\Lambda$.  The {\it Scarf complex} of $I$ is the labelled
subcomplex $\Delta_I \subseteq X$ determined by
$$
  \Delta_I \ =\ \{F \in X \mid \aa_F = \aa_G \Rightarrow F = G\}.
$$
The free complex $\FF_{\Delta_I}$ it supports is acyclic as well as
minimal, and coincides with the hull resolution of $I$ by \cite{BS},
Theorem~2.9.  See also Example~\ref{fig:bounded}. \hfill $\Box$
\end{example}

\begin{defn}[The cohull resolution] \label{defn:cohull}
The {\rm cohull resolution} ${\rm cohull}_\aa(I)$ of an ideal $I$ with
respect to $\aa \succeq \aa_I$ is defined to be the relative cocellular
resolution $\FF^{\,(X,X_U)}[-\aa -{\bf 1}](1-n)$, where $X$ is the hull
complex of $I^\aa + \mm^{\aa + {\bf 1}}$.  The {\rm canonical cohull
resolution}, or simply {\rm the cohull resolution} ${\rm cohull(I)}$ of
$I$ is obtained by taking $\aa = \aa_I$.
\end{defn}
The cohull resolution, like the hull resolution, is a possibly nonminimal
resolution that preserves some of the symmetry (in the generators and
irreducible components) of an ideal.

There are some geometric properties of hull resolutions of artinian
ideals that make cohull resolutions a little more tangible.  Suppose, for
instance, that $J$ is an artinian monomial ideal, with $x_1^{d_1},
\ldots, x_n^{d_n}$ among its minimal generators.  Choose $t > (n+1)!$,
and let $v_1, \ldots, v_n$ be the vertices of the polyhedron $P_t$
determined by these minimal generators.  The vertices $\{v_i\}$ of $P_t$
span an affine hyperplane which will be denoted by $H$.

Fix a strictly positive inner normal $\varphi_G$ for each $G \in {\rm
hull}(J)$.  Recall that $P_t$ is contained in the (closed) polyhedron
${\bf 1} + \RR_+^n$ (since monomials in $S$ have no negative exponents).
Each face $G \in$ hull$(J)$ spans an affine space which does not contain
the vector ${\bf 1} \in \RR^n$ because the hyperplane containing $G$ and
normal to $\varphi_G$ does not contain ${\bf 1}$.  Therefore the
projection $\pi$ from the point ${\bf 1}$ to the hyperplane $H$ induces a
homeomorphism hull$(J) \to \pi($hull$(J))$.  In fact,

\begin{prop} \label{prop:subdivision}
If $J$ is artinian, $\pi({\rm hull}(J))$ is a regular polytopal
subdivision of the simplex $H \cap P_t$.
\end{prop}
{\it Proof:\ } That $H \cap P_t \subset {\bf 1} + \RR_+^n$ is a simplex
follows because it is convex and contains $v_1,\ldots,v_n$.  Now $\pi$
induces a map of the boundary $\partial P_t \to H \cap P_t$ which is
obviously surjective.  Suppose that $\pi(\ww)$ is in the interior of $H
\cap P_t$ for some $\ww \in \partial P_t$.  It is enough to show that if
a nonzero support functional $\varphi$ attains its minimum on $P_t$ at
$\ww$ then $\varphi$ is strictly positive.  All coordinates of $\varphi$
are $\geq 0$ {\it a priori} because it attains a minimum on $P_t$; but if
the $i^{\rm \,th}$ coordinate of $\varphi$ is zero then $\<\varphi,v_i\>
< \<\varphi, \ww\>$ and $\varphi$ cannot be minimized at $\ww$.
\hfill
$\Box$

\begin{remark} \rm
This generalizes the result \cite{BPS}, Corollary~5.5 for generic
artinian monomial ideals, in view of \cite{BS}, Theorem~2.9.  Regular
subdivisions here are as in \cite{Zie}, Definition~5.3.
\end{remark}

We arrive at the following characterization for artinian hull complexes:
\begin{thm} \label{thm:boundary}
If $X$ is the hull complex of an artinian monomial ideal, then $|X|$ is a
simplex and the negatively unbounded complex $X_U$ is the topological
boundary of $X$.
\end{thm}
{\it Proof:\ } By the previous proposition, it suffices to show that a
face $G$ of the hull complex of any (not necessarily artinian) ideal has
a label without full support if and only if it is contained in the
topological boundary of the shifted positive orthant ${\bf 1} + \RR_+^n$.
But this holds because the $i^{\,\rm th}$ coordinate of $\aa_G$ is zero
if and only if every vertex of $G$ (and hence every point in $G$) has
$i^{\,\rm th}$ coordinate $1$.
\hfill
$\Box$
\vskip 2mm

Although cohull resolutions are relative cocellular by definition, they
can frequently be viewed as cellular resolutions, as well.  In fact, with
a slight weakening of the notion of labelled cell complex, all cohull
resolutions are weakly cellular.  To be precise, define a {\it weakly
labelled cell complex} to be the same as a labelled cell complex, except
that instead of requiring that the label $\aa_F$ be equal to the join
$\bigvee_{v \in F} \aa_v$, we require only that $\aa_F \preceq \bigvee_{v
\in F} \aa_v$ whenever $\dim F > 0$.  A free complex or resolution is
called {\it weakly cellular} if it is supported on a weakly labelled cell
complex.

\begin{thm} \label{thm:weakly}
The cohull resolution of $I$ with respect to $\aa$ is weakly cellular for
any $\aa \succeq \aa_I$.
\end{thm}
{\it Proof:\ } Let $J = I + \mm^{\aa + {\bf 1}}$ and assume the notation
from after Definition~\ref{defn:cohull}.  Define $Q_t$ to be the
intersection of $P_t$ with the closed half-space containing the origin
and determined by the hyperplane $H$.  Then $Q_t$ is a polytope which may
also be described as the convex hull of (all of) the vertices of $P_t$.
Furthermore, the bounded faces of $P_t$ are simply those faces of $Q_t$
which admit a strictly positive inner normal.  Thus $X := {\rm hull}(J)$
is a subcomplex of the boundary complex of $Q_t$, as is the boundary
$\partial X$.

Let $Y \subset \partial Q_t$ be the subcomplex generated by the facets of
$Q_t$ whose inner normal is {\it not} strictly positive.  Denote chain
and relative cochain complexes over $k$ by ${\cal C}\ldot(-)$ and ${\cal
C}\hdot(-,-)$.  Then $Y \cap X = \partial X$ and the ${\cal
C}\hdot(Q_t,Y) = {\cal C}\hdot(X,\partial X)$.  For elementary reasons,
${\cal C}\hdot(Q_t,Y) \cong {\cal C}\ldot(X^\vee)$ for some subcomplex
$X^\vee$ of the polar polytope $Q_t^\vee$ (use, for instance, the methods
of \cite{Zie}, Sections~2.2--2.3).  Note that the isomorphisms will exist
regardless of the incidence functions in question, by \cite{BH},
Theorem~6.2.2.  That $X^\vee$ is weakly labelled follows from the
isomorphism ${\cal C}\ldot(X^\vee) \cong {\cal C}\hdot(X,\partial X)$ and
the remark after Definition~\ref{defn:relcocell}.  Indeed, the condition
$F \supseteq G \Rightarrow \aa_F \succeq \aa_G$ for faces $F,G \in
(X,\partial X)$ is equivalent to the condition $F^\vee \subseteq G^\vee
\Rightarrow -\aa_F \preceq -\aa_G$ for faces of $X^\vee$, and this need
only be applied when $F$ is a facet containing $G$ and $F^\vee$ is a
vertex of $G^\vee$.
\hfill
$\Box$

\begin{prop} \label{prop:minimal}
If a weakly cellular resolution is minimal, it is cellular.  In
particular, if a cohull resolution is minimal, it is cellular.
\end{prop}
{\it Proof:\ } Let $(\widetilde{\FF},\widetilde{\partial})$ denote the
augmented complex $\FF_X \to I \to 0$, where $X$ is a weakly labelled
complex supporting a free resolution of $I$.  We show that if $G \in X$
then $\aa_G \succ \bigvee_{v \in G} \aa_v$ implies $\FF$ is not minimal.
This is vacuous if $\dim G = 0$, so assume $\dim G$ has minimal dimension
$\geq 1$, and suppose that $\aa_G - \ee_i \succeq \bigvee_{v \in G}
\aa_v$.  Then $\widetilde{\partial}(G) = x_iy$ for some $y \in
\widetilde{\FF}$ because $\dim G$ is minimal.  It follows that
$x_i\widetilde{\partial}(y) = \widetilde{\partial}(x_iy) = 0$, whence
$\widetilde{\partial}(y) = 0$ because $\widetilde{\FF}$ is torsion-free.
Thus $\widetilde{\partial}(G) \in x_i{\rm ker}(\widetilde{\partial})
\subseteq \mm \cdot {\rm ker}(\widetilde {\partial})$ does not represent
a minimal generator of ${\rm ker}(\widetilde{\partial})$ by Nakayama's
Lemma for graded modules.
\hfill
$\Box$
\vskip 2mm

\noindent
\begin{remark} \rm
For cohull resolutions the proposition is probably true without the
hypothesis of minimality, but a proof (which would likely be geometric
instead of algebraic) has not been found.  In particular, all cohull
resolutions in the examples below are cellular.  Cellularity of the
cohull resolution is equivalent to the following more concrete statement:
the label on any interior face of the hull complex of an artinian ideal
is the greatest common divisor of the labels on the facets that contain
it.
\end{remark}

\begin{example} \label{ex:permut3} \rm
{\sl (continuation of Example~\ref{ex:permut})\ } The minimal resolution
of the permutahedron ideal $I$ of Example~\ref{ex:permut} is, by
\cite{BS}, Example~1.9, the hull resolution, which is supported on a
permutahedron.  The minimal resolution of $I + \mm^{(n+1){\bf 1}}$ is
also the hull resolution, and is supported on the complex $X$ which may
described as follows.

There are two kinds of faces of $X$.  The first kind are those that make
up the boundary $\partial X$; these are indexed by the proper nonempty
faces $F \in \Delta$ and have vertices $t^{(n+1) \ee_i} \in P_t$ for $i
\in F$ (recall from Section~\ref{section:definitions} that $\ee_i$
denotes the $i^{\,\rm th}$ basis vector of $\ZZ^n$ and $\Delta =
\{1,\ldots,n\}$ is the $(n\!-\!1)$-simplex).  On the other hand, the
interior \hbox{$p$-faces} of $X$ are in bijection with the chains
\begin{equation} \label{eqn:face}
  \emptyset \,\prec\, F_1 \,\prec\, F_2 \,\prec\, \cdots \,\prec\,
  F_{n-p}
\end{equation}
of faces of $\Delta$, where $F_{n-p}$ might (or might not) equal
$\Delta$.  Note that the interior faces of $X$ for which $F_{n-p} =
\Delta$ are faces of the permutahedron itself.

More generally, an interior \hbox{$p$-face} $G$ given by (\ref{eqn:face})
for which $F_{n-p} \neq \Delta$ is affinely spanned by the permutahedral
\hbox{$(p\!-\!1)$-face} $G': \emptyset \prec F_1 \prec \cdots \prec
F_{n-p} \prec \Delta$ and the ``artinian'' vertices $\{t^{(n+1)\ee_i}
\mid i \not\in F_{n-p}\}$ of $P_t$.  In fact, a functional which attains
its minimum (in $P_t$) on $G$ may be produced directly.  For this
purpose, define for any $F \in \Delta$ the functional $F^\dagger$ on
$\RR^n$ to be the transpose of $F$; i.e.\ $\<F^\dagger,\ee_i\> = 1$ if $i
\in F$ and zero otherwise.  Then the functional $\varphi_\epsilon := {\bf
1}^\dagger + \epsilon \sum_{j=1}^{n-p} F_j^\dagger$ attains its minimum
(in $P_t$) on $G'$ for all $0 < \epsilon <\!\!< 1$.  But for $\epsilon
>\!\!> 0$ we have $\<\varphi_\epsilon, t^{(n+1)\ee_i}\> <
\<\varphi_\epsilon, G'\>$ whenever $i \not\in F_{n-p}$.  Thus we can
choose the unique $\epsilon$ that makes $\<\varphi_\epsilon,
t^{(n+1)\ee_i}\> = \<\varphi_\epsilon, G'\>$ for all $i \not\in F_{n-p}$,
so that $\varphi_\epsilon$ attains its minimum on $G$.

It is easy to check that the labels on the faces of $X$ are distinct,
whence $\FF_X$ is the minimal resolution of $I + \mm^{(n+1){\bf 1}}$ by
\cite{BS}, Remark~1.4.  In particular, the irredundant irreducible
components of $I$ are in bijection with the facets of $X$ by
Theorem~\ref{thm:irr}, and the generators of the tree ideal $I^\vee$
are given by $x^{(n+1) {\bf 1} - \aa_G}$ for facets $G \in X$.  This
recovers the generators for $I^\vee$ in Example~\ref{ex:permut}.

Retaining earlier notation, the face $G$ has dimension $1 + \dim(G')$.
Thus the \hbox{$p$-faces} of $X$ are in bijection with the collection of
\hbox{$p$-\ and $(p\!-\!1)$-faces} of the permutahedron.  In fact, the
(unlabelled) pair $(|X|,|\partial X|)$ has the same faces as the pair
$\Bigl( \partial(v * Y),v \Bigr)$ consisting of the boundary of the cone
over the permutahedron $Y$ rel the apex of the cone.
\begin{figure}
\vbox{\vskip 1mm \baselineskip 0pt 
\begin{centering}\includegraphics{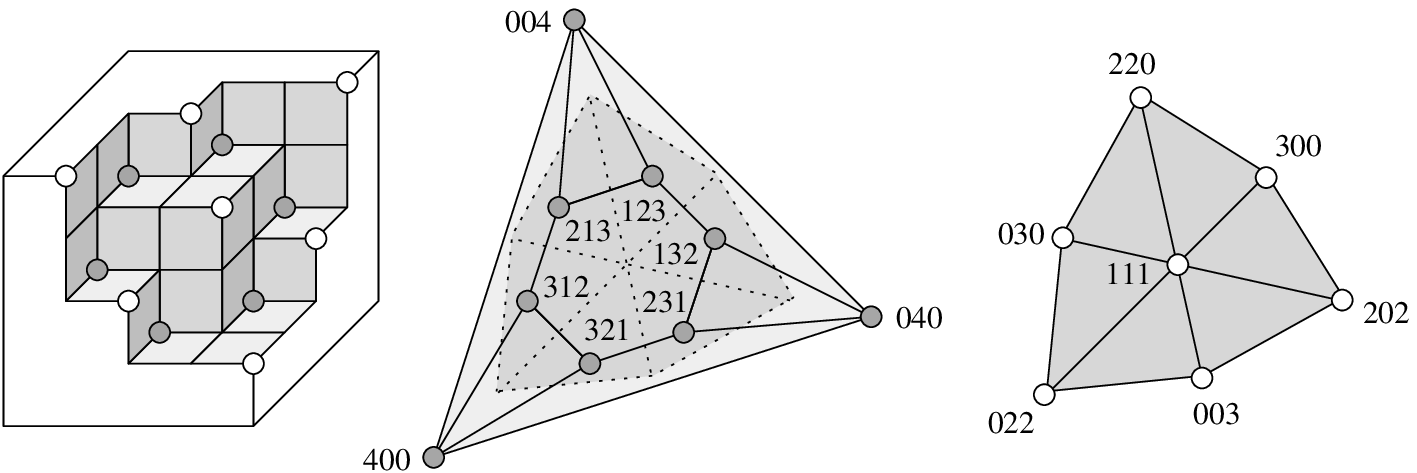}\\
\hfill\hfill $I$ \hfill\hfill\hfill $X =$ hull$(I + \mm^{(4,4,4)})$
\hfill\hfill $X^\vee = $ cohull$(I^\vee)$ \ \ \hfill
\end{centering}
\vskip 2mm
\noindent
Figure~6: $I$ and $I^\vee$ are the permutahedron and tree ideals when
$n = 3$.  The complex $X$ is the (labelled) regular polytopal subdivision
of the simplex promised by Proposition~\ref{prop:subdivision}.  Overlayed
on this figure is the dual complex $X^\vee$ (without its labelling).  At
right, $X^\vee$ is shown with its labelling, which is $\ZZ^n$-shifted as
per Theorem~\ref{thm:relcocell}.  Turn the picture over for the staircase
of $I^\vee$.
}
\end{figure}
The cellular complex $X^\vee$ supporting the cohull resolution of the
tree ideal $I^\vee$ is therefore easy to describe.  Let $Y$ be the
permutahedron in $\RR^n$ and $Y^\vee$ its polar.  Then $X^\vee$ is the
cone over $\partial Y^\vee$ from the barycenter of $Y^\vee$.  The
vertices $G^\vee$ of $X^\vee$, which are labelled by the generators of
$I^\vee$, almost all correspond to the facets $G'$ of $Y$ (whose
labellings are as above).  Only the apex of the cone is an exception,
corresponding instead to the interior of $Y$.  The case $n=3$ is depicted
in Figure~6; it should be noted that the equality $Y = Y^\vee$ is only
because $Y$ is 2-dimensional, not some more general self-duality.

Now cohull$(I^\vee)$ is a cellular resolution of $I^\vee = I^\vee +
\mm^{\aa_I + {\bf 1}}$, so we can dualize this cellular resolution using
Theorem~\ref{thm:relcocell} again.  This yields a minimal relative
cocellular resolution of $I$, which is seen to be cellular and
(coincidentally?) equal to hull$(I)$.
\hfill
$\Box$
\end{example}

Recall from Section~\ref{section:definitions} that an ideal is {\it
cogeneric} if it is Alexander dual to a generic ideal.  The minimal
resolution of such an ideal was introduced in \cite{Stu}, where it was
dubbed the {\it co-Scarf resolution}.  The next theorem, along with the
proof of Theorem~\ref{thm:weakly} above, explains why the construction in
\cite{Stu} involved a subcomplex of the boundary of the simple polytope
dual to the simplicial polytope of which the Scarf complex is a
subcomplex.  The theorem is a direct consequence of
Theorem~\ref{thm:relcocell}, Example~\ref{ex:scarf}, and
Proposition~\ref{prop:minimal}.
\begin{thm} \label{thm:coscarf}
Any cohull resolution of a cogeneric monomial ideal is minimal and
cellular.
\end{thm}
\begin{remark} \rm
That the co-Scarf resolution is cellular as opposed to weakly cellular
was assumed in \cite{BS}, Example~1.8 but overlooked in \cite{Stu}.
\end{remark}
\begin{figure}[!hbp]
\mbox{}
\hfill hull$(I^\aa)$
\hfill \hfill \hfill \hfill $I^\aa$
\hfill \hfill \hfill hull$(I^\aa + \mm^{\aa + {\bf 1}})$
\hfill cohull$_\aa(I)$
\hfill 
\hfill \hfill \hfill \hfill \mbox{}

\vskip 1mm
\noindent
\hfill \includegraphics{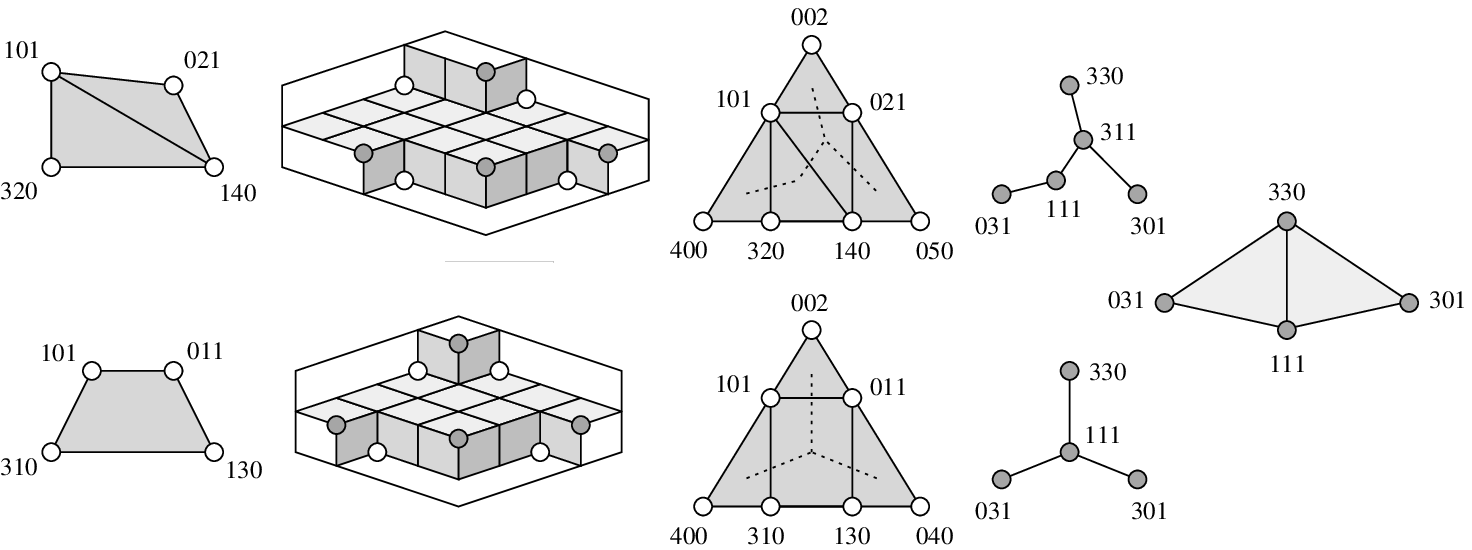} \hfill \mbox{}

\vskip 1mm
\mbox{}
\hfill hull$((I^\aa)^\vee{}^\vee)$
\hfill \hfill \hfill \hfill \hfill \hfill $I^\vee$
\hfill \hfill \hfill \hfill \hfill hull$(I^\vee + \mm^{\aa_I + {\bf 1}})$
\hfill \hbox{\hskip 1ex cohull$(I)$}
\hfill \hfill 
\raisebox{32ex}[0mm][0mm]{\makebox[0mm][l]{$\phantom{_I\,}\aa = (3,4,1)$}}
\raisebox{6ex}[0mm][0mm]{\makebox[0mm][l]{$\aa_I = (3,3,1)$}}
\raisebox{11ex}[0mm][0mm]{\hbox{\hskip -1ex hull$(I)$ \hskip 1ex}}
\hfill \hfill \hfill \mbox{}
\vskip 2mm
\noindent
{\hfill Figure~7 \hfill \mbox{}}
\vskip 5ex
\vbox{
\noindent
\mbox{}
\hbox{\hskip -1ex $X$, minimally resolves $I^\vee + \mm^{(5,4,3)}$}
\hfill \hbox{\hskip 2.5ex hull$(I)$ \hskip 1ex}
\hfill \hfill \hfill \hfill \hbox{\hskip 2ex hull$(I^\vee \! +
					\mm^{(5,4,3)})$ \hskip -2ex}
\hfill \hfill \hfill \mbox{}

\vskip 1mm
\noindent
\mbox{} \hfill \includegraphics{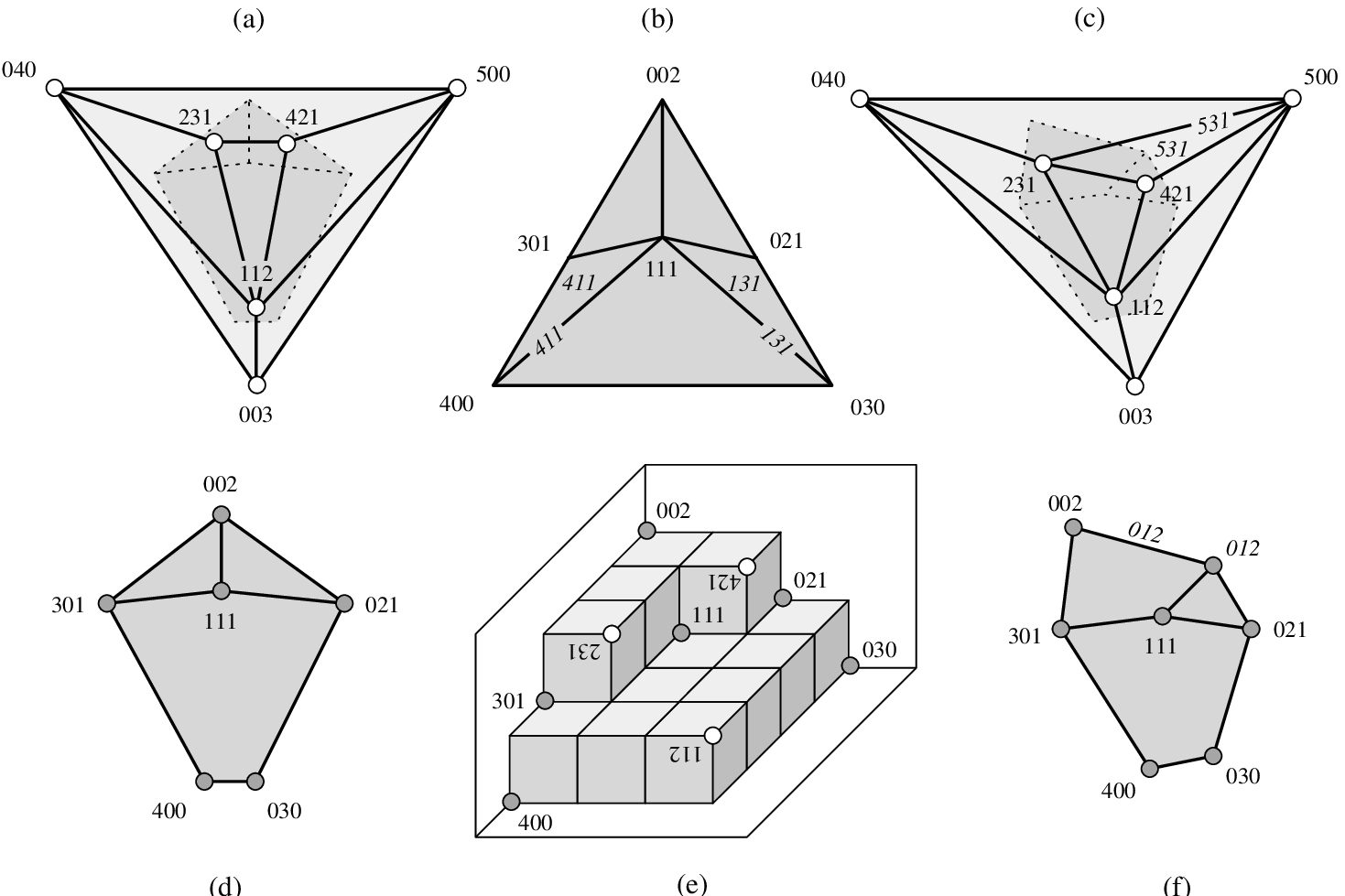} \hfill \mbox{}

\vbox{\baselineskip 0pt
\vskip 1ex
\noindent
\mbox{}
\hfill \hbox{$Y$, minimally resolves $I$}
\hfill \hbox{\,staircase of $\:I\phantom{^\vee}\!\!$ rightside up}
\hfill \hfill \hfill \hbox{\hskip 3ex cohull$(I)$ \hskip -3ex}
\hfill \hfill \hfill \mbox{}
\newline
\mbox{}
\hfill \phantom{\hbox{$Y$, minimally resolves $I$}}
\hfill \hbox{$\;$staircase of $I^\vee$ upside down}
\hfill \hfill \hfill \phantom{\hbox{\hskip 3ex cohull$(I)$ \hskip -3ex}}
\hfill \hfill \hfill \mbox{}
}
\vskip 2mm
\noindent
{\hfill Figure~8 \hfill}
\vskip 2mm}
\end{figure}
\begin{example} \rm
It is possible for the hull and cohull resolutions to coincide for a
given ideal $I$.  For instance, this occurs if $I = \mm$; or if $I$ is
simultaneously generic and cogeneric (which turns out to be pretty hard
to accomplish!); or if $I$ is the permutahedron ideal in 3 variables.
Conjecturally, the hull and cohull resolutions should coincide for
permutahedron ideals of all dimensions.  \hfill $\Box$
\end{example}
\begin{example} \rm
Of course, it is also possible for the hull and cohull resolutions to be
very different.  For instance, the cohull resolution of the ideal
$I^\vee$ from Examples~\ref{ex:alexdual} and~\ref{ex:coscarf} is the
co-Scarf resolution, which is cellular and supported on an octagon with
only one maximal face (dualize the picture in Figure~5).  On the other
hand, ${\rm hull}(I^\vee)$ is a triangulation of the same octagon.
\hfill $\Box$
\end{example}
\begin{example} \label{ex:tighten} \rm
The canonical cohull resolution can differ from a noncanonical cohull
resolution: let $I = \<x^3z,xyz,y^3z,x^3y^3\>$ and $\aa = (3,4,1)$, so
$I^\aa = \<xz, x^3y^2, xy^4, y^2z\>$ and $I^\vee = \<xz, x^3y, xy^3,
yz\>$.  Since hull$(I)$ is not minimal, we look elsewhere for the minimal
resolution of $I$.  But hull$(I^\aa + \mm^{\aa + {\bf 1}})$ is not
minimal, and the failure of minimality occurs in such a way that
cohull$_\aa(I)$ is also not minimal.  On the other hand, the offending
nonminimal edge is not present in hull$(I^\vee + \mm^{\aa_I + {\bf 1}})$,
and this resolution is minimal.  It follows that cohull$(I)$ is minimal.
Note how the passage from $I^\aa$ to $(I^\aa)^\vee{}^\vee = I^\vee$
``tightens'' the hull resolution of $I^\aa$ to make the nonminimal edge
disappear in hull$((I^\aa)^\vee{}^\vee)$, cf.\ the remark after
Corollary~\ref{cor:Iaa=I}.

The labelled complexes supporting these resolutions are all depicted in
Figure~7, where the resolutions with black vertices are drawn ``upside
down'' to make their superimposition on the staircase diagram for $I$
easier to visualize.  Observe that a staircase diagram for $I$ can be
obtained by turning over the staircase diagram for either $I^\aa$ or
$I^\vee$, although these result in different ``bounding boxes'' for $I$.
Note that all of the complexes, particularly the cohull complexes, are
labelled and not just weakly labelled.  \hfill $\Box$
\end{example}
\begin{example} \rm
Finally, an example to illustrate that not all cellular resolutions come
directly from hull and cohull resolutions, so that the algebraic
techniques to prove exactness in Section~\ref{section:limits} prove a
stronger duality for resolutions than a geometric treatment such as that
in \cite{BS} or \cite{Stu} could provide.  All of the labelled cellular
complexes from this example are depicted in Figure~8.

Let $I = \<z^2,x^3z,x^4,y^3,y^2z,xyz\>$, so that $I^\vee =
\<xyz^2,x^2y^3z,x^4y^2z\>$.  Then hull$(I)$ and cohull$(I)$ are not
minimal (the offending cells have italic labels); moreover, ${\rm
cohull}_\aa(I) = {\rm cohull}(I)$ for all $\aa \succeq \aa_I = (4,3,2)$.
Nonetheless, the minimal resolution $\FF_X$ of $I^\vee + \mm^{(5,4,3)}$
is cellular, so Theorem~\ref{thm:relcocell} applies, yielding a minimal
relative cocellular resolution $\FF^{\,(X,X_U)}[-(5,4,3)](-2)$ for $I$.
In fact, this relative cocellular resolution is cellular, supported on the
labelled cell complex $Y$.  \hfill $\Box$
\end{example}

\end{section}
\begin{section}{Deformations and limits of resolutions}

\label{section:limits}

The final item on the agenda is the proof of Theorem~\ref{thm:relcocell}.
To that end, the goal of this section is Theorem~\ref{thm:coresolution},
which is actually a little more general than Theorem~\ref{thm:relcocell}.
It can be viewed as the result of applying a limiting process to a
collection of pairs of linked artinian monomial ideals that are
deformations of a given pair.  The entire section is a setup to apply a
limit to Proposition~\ref{prop:ext}, and is another manifestation of the
kinship of Alexander duality and other types of duality for Gorenstein
rings.  The maps $f_\bb$ in the following definition accomplish the
deformations.

\begin{defn} \label{defn:fb}
Define the map $f_\bb \colon \ZZ^n \to \ZZ^n$ for $\bb \succeq {\bf 0}$
by the coordinatewise formula
$$f_\bb(\cc)_i\ =\
\left\{
	\begin{array}{ll}
	c_i - b_i & {\rm if}\ c_i \leq 0 \cr
	c_i	  & {\rm if}\ c_i \geq 1
	\end{array}
\right.
$$
To avoid messy exponents we also let $f_\bb(x^\cc) = x^{f_\bb(\cc)}$.
Whenever the symbol $f_\bb$ is written, it will be assumed that $\bb
\succeq {\bf 0}$.
\end{defn}

\begin{prop} \label{prop:fbI=SbI}
If $I \subseteq S$ is any monomial ideal then $\<f_\bb(I)\> = S[\bb] \cap
\widetilde{I}$.
\end{prop}
{\it Proof:\ } It is clear from the definition that $f_\bb(\cc) \succeq
-\bb$ if $\cc \succeq {\bf 0}$, whence $\<f_\bb(I)\> \subseteq S[\bb]$.
Since also $f_\bb(\cc)^+ = f_\bb(\cc^+)$, we conclude that $\<f_\bb(I)\>
\subseteq \widetilde{I}$ as well.  For the reverse inclusion, assume
$x^\cc \in S[\bb] \cap \widetilde{I}$.  Then $f_\bb(x^{\cc^+}) \in
f_\bb(I)$ and divides $x^\cc$ because $f_\bb(\cc^+) \preceq \cc$ whenever
$\cc \succeq -\bb$, a fact which is easily seen from the definition.
\hfill
$\Box$ 
\vskip .2cm

Recall from Section~\ref{section:cellular} that the join of $\cc, \cc'
\in \ZZ^n$ is the componentwise maximum.
\begin{lemma} \label{lemma:joins}
The map $f_\bb$ preserves joins; that is, $f_\bb(\cc \vee \cc') =
f_\bb(\cc) \vee f_\bb(\cc')$. \hfill $\Box$
\end{lemma}
This lemma is important because of the next proposition, originally due
to D.\ Bayer.  Let $X$ be a labelled cell complex, and suppose $f \colon
\ZZ^n \to \ZZ^n$ is a map respecting joins.  Denote by $f(X)$ the
labelled cell complex which is obtained by applying $f$ to the labels on
the faces of $X$.  Thus $G \in f(X)$ is labelled by $f(\aa_G)$ whenever
$G \in X$ is labelled by $\aa_G$.
\begin{prop} \label{prop:joins}
Let $\FF_X$ be a cellular resolution of a finitely generated module $M
\subseteq T$.  If $f \colon \ZZ^n \to \ZZ^n$ preserves joins then
$\FF_{f(X)}$ is a resolution of $\<f(M)\>$.
\end{prop}
{\it Proof:\ } Note that because $f$ respects joins the effect of $f$ is
determined by its effect on the vertex labels.  Similarly, $\<f(M)\> =
\<f(x^\bb) \mid \bb$ is a vertex label of $X\>$.  Thus one only needs to
check that $\FF_{f(X)}$ is acyclic.  It suffices to check that $X_B(\bb)$
is acyclic for all $\bb \in \ZZ^n$, by the acyclicity criterion of
\cite{BS}, Proposition~1.2.

Suppose, then, that $\alpha$ is a cycle of the reduced chain complex of
$|f(X)_B(\bb)|$.  Then $\alpha$ also represents a cycle of $|X|$.  Let
$\cc$ be the join of the labels on the faces in the support of $\alpha$,
considered as faces of $X$.  Since $f$ preserves joins, $f(\cc) \preceq
\bb$ and $|X_B(\cc)| \subseteq |f(X)_B(\bb)|$.  Now $\alpha$ is a
boundary in the reduced chain complex of $|X_B(\cc)|$ by \cite{BS},
Proposition~1.2, and it follows that $\alpha$ is also a boundary in the
reduced chain complex of $|X_B(\bb)|$, completing the proof.
\hfill
$\Box$
\begin{cor} \label{cor:fbX=SbI}
If $\FF_X$ is a cellular resolution of $I$ then $\FF_{f_{ \scriptstyle
\bb }(X)}$ is a cellular resolution of $S[\bb] \cap \widetilde{I}$.
\hfill $\Box$
\end{cor}

Keeping the notation of the corollary we can augment the complex
$\FF_{f_{\scriptstyle \bb}(X)}$ to a resolution of $S[\bb]/S[\bb] \cap
\widetilde{I}$, homologically shifted down 1, by adding a summand
$S[\bb]$ in homological degree $-1$.  We denote this augmented resolution
by $\FF_\bb(X)$, and we let $\FF^{\,\bb}(X) := \FF_\bb(X)^\Ast$, with
differential $\delta_\bb$.  The generator of the summand
$S[-f_\bb(\aa_F)] \subseteq \FF_\bb(X)$ corresponding to the face $F$
will be denoted by $F_\bb$, while the generator of $S[f_\bb(\aa_F)] =
S[-f_\bb(\aa_F)]^\Ast \subseteq \FF^{\,\bb}(X)$ will be denoted by
$F^\bb$.  Keep in mind that $F^\bb$ is in $\ZZ^n$-graded degree
$-f_\bb(\aa_F)$.

We will soon be defining maps between the $\FF^{\,\bb}(X)$ for various
$\bb$, and the following lemma, particularly part~(ii), will be the tool
used to prove that these maps are well-defined, commute with the
differentials, and form an inverse system.
\begin{lemma} \label{lemma:f-properties}
If $\bb \succeq \bb' \succeq {\bf 0} $ then 
$$
\begin{array}{ll}
(i)
&	f_\bb = f_{\bb - \bb'} \,{\scriptstyle \circ}\, f_{\bb'}\,,
\cr
(ii)
&	f_{\bb'}(\cc) - f_\bb(\cc) = \cc - f_{\bb-\bb'}(\cc)\,.
\end{array}
$$
\end{lemma}
{\it Proof:\ } Plug and chug, using the equality $f_\bb(\cc)^+ = \cc^+$
for (i).
\hfill
$\Box$ 
\vskip .2cm

\begin{lemma} \label{lemma:injection}
For every $\bb \succeq \bb' \succeq {\bf 0} $ we have an injection
$\varphi_{\bb,\bb'} \colon \FF^{\,\bb}(X) \hookrightarrow
\FF^{\,\bb'}(X)$ of chain complexes sending $F^\bb$ to ${\displaystyle
\phantom{(} {\textstyle m_F} \phantom{)} \over {\textstyle f_{\bb -
\bb'}(m_F)} }F^{\bb'}$.
\end{lemma}
{\it Proof:\ } There are two aspects to the proof: (i) the given map is
an injection of homologically graded modules which (as a map of
$\ZZ^n$-graded modules) has degree ${\bf 0}$, and (ii) the injections
commute with the differentials.  The first follows from the equality
$-f_\bb(\aa_F) = -f_{\bb'}(\aa_F) + \aa_F - f_{\bb-\bb'}(\aa_F)$ which is
easily seen to be equivalent to Lemma~\ref{lemma:f-properties}(ii) when
$\cc = \aa_F$.  The second is a longer calculation directly from the
definition of the differentials $\delta_\bb$ and $\delta_{\bb'}$ of the
chain complexes $\FF_\bb$ and $\FF_{\bb'}$.  

The definitions imply that $\delta_\bb$ is just the transpose of the
differential from the cellular free complex as defined in \cite{BS}.
Thus, $\delta_\bb(F^\bb) = \sum_{G \in X} \varepsilon(G,F) {f_\bb(m_G)
\over f_\bb(m_F)} G^\bb$, where $\varepsilon$ is the incidence function
defining the differential of $X$.  Note that $\varepsilon(G,F)$ is
nonzero only if $G \supseteq F$.  We have
\begin{eqnarray*}
	\delta_{\bb'} \circ \varphi_{\bb,\bb'}(F^\bb)	&
	=	& 
	\displaystyle
	\sum_{G \in X} \varepsilon(G,F) { {\textstyle f_{\bb'}(m_G)}
	\over {\textstyle f_{\bb'}(m_F)} } \cdot { {\textstyle m_F}
	\over {\textstyle f_{\bb - \bb'}(m_F)} }\,G^{\,\bb'}
\cr
\cr
		& 
	=	& 
	\displaystyle
	\sum_{G \in X} \varepsilon(G,F) { {\textstyle m_G} \over
	{\textstyle f_{\bb - \bb'}(m_G)} } \cdot { {\textstyle
	f_{\bb}(m_G)} \over {\textstyle f_{\bb}(m_{F})} }\,G^{\,\bb'}
\cr
\cr
		&
	=	& 
	\varphi_{\bb,\bb'} \circ \delta_\bb(F^\bb),
\end{eqnarray*}
where the transition from the first line to the second is accomplished by
two applications of Lemma~\ref{lemma:f-properties}(ii). 
\hfill
$\Box$ 
\vskip .2cm

\begin{lemma}
If $\bb \succeq \bb' \succeq \bb'' \succeq {\bf 0}$ then
$\varphi_{\bb,\bb''} = \varphi_{\bb',\bb''}\,{\scriptstyle \circ}\,
\varphi_{\bb,\bb'}$\,.
\end{lemma}
\noindent
{\it Proof:\ } We need only check the equality as maps of modules.  The
proof again uses property (ii) from Lemma~\ref{lemma:f-properties}, and
it involves manipulations similar to those in the proof of
Lemma~\ref{lemma:injection}.
\hfill 
$\Box$
\vskip .2cm

These lemmata show that we have an inverse system of complexes of free
modules, so it is natural now to take the inverse limit.  With
$\FF^{\,t}(X) := \FF^{\,t \cdot {\bf 1}}(X)$ we can simplify a little
since the inverse systems $\{\FF^{\,\bb}(X)\}_{\bb \succeq 0}$ and
$\{\FF^{\,t}(X)\}_{t \in \NN}$ are cofinal, so that their limits are the
same.  We take this opportunity to note that our inverse limits, when
taken in the category of $\ZZ^n$-graded objects and degree zero maps,
will be denoted by $^\Ast\invlim{}\!$, and that $S$ is complete in this
category.  Recall that, for our inverse system $\{\FF^{\,t}(X)\}_{t \in
\NN}$ of chain complexes, for instance, this is defined as
$$
	{}^\Ast\Invlim{t} \FF^{\,t}(X)\ \, =\ 
	\bigoplus_{\cc \in \ZZ^n}\; \Invlim{t} \FF^{\,t}(X)_\cc\,,
$$
where the inverse limits on the right are in the category of chain
complexes of $k$-vector spaces.

At each stage in the inverse system, $f_\bb$ moves the labels on $X_U$
away from the first orthant, in negative directions, turning any zeros
into arbitrarily large negative integers (hence the name ``negatively
unbounded'' for the subcomplex $X_U$ of $X$).  Then $S$-duality makes the
negative integers positive.  Thus the maps $f_\bb$, combined as they are
with $S$-duality in the definition of $\FF^\bb$, create irreducible
components of $I^\aa$ from those generators of $I$ which do not have full
support by pushing the zeros out to (positive) infinity.  In the limit,
the vertices defining those generators disappear.  This provides the
intuition for the next result.
\begin{thm} \label{thm:limit}
$\FF^{\,(X,X_U)} = {}^\Ast\invlim{t} \FF^{\,t}(X)$.
\end{thm}
{\it Proof:\ } The first observations are that $\FF^{\,(X,X_U)}$ is a
subcomplex of $\FF^{\,t}$ for all $t$, and that the maps
$\varphi_{t\,,\,t'} := \varphi_{t \cdot {\bf 1}\,,\, t' \cdot {\bf 1}}$
defining the inverse system restrict to the identity on
$\FF^{\,(X,X_U)}$.  This is because of the way $f_t := f_{t \cdot {\bf
1}}$ is defined:
\begin{equation} \label{eqn:fixed-monomials}
f_{t-t^{\,\prime}}(m_F) = m_F\ \ \Longleftrightarrow\ \ t = t^{\,\prime}\
{\rm or}\ F \not \in X_U
\end{equation}
because $f_\bb(\cc)_i = c_i$ for all $i$ precisely when $\cc \succeq {\bf
1}$.  Thus we have, for all $t \geq 0$, 
\begin{equation} \label{ses:invlims}
0 \,\to\, \FF^{\,(X,X_U)} \,\to\, \FF^{\,t}(X) \,\to\, \FF^{\,t}(X_U)
\,\to\, 0
\end{equation}
giving rise to a corresponding exact sequence of inverse systems.  To be
more precise, the maps $\{\varphi_{t\,,\,t^{\,\prime}}\}$ from the
inverse system $\{\FF^{\,t}(X)\}_{t \in \NN}$ induce maps
$\{\psi_{t\,,\,t^{\,\prime}} \colon \FF^{\,t}(X_U) \to
\FF^{\,t^{\,\prime}}(X_U)\}_{t \geq t^{\,\prime}}$ which make
$\{\FF^{\,t}(X_U)\}_{t \in \NN}\,$ into an inverse system.

It is readily seen that the maps $\psi_{t\,,\,t^{\,\prime}}$ in the
inverse system are injections, so the limit is ${}^\Ast\invlim{t}
\FF^{\,t}(X_U) = \bigcap_{\,t}\,\psi_{t\,,\,0} \bigl( \FF^{\,t}(X_U)
\bigr)$.  Furthermore, (\ref{eqn:fixed-monomials}) implies that
$\psi_{t\,,\,t^{\,\prime}} \bigl( \FF^{\,t}(X_U) \bigr) \subseteq
\mm\FF^{\,t^{\,\prime}}(X_U)$ if $t > t^{\,\prime}$.  It follows from the
Krull intersection theorem that the inverse limit is zero.  Since the
inverse limit is always left exact our exact sequence of inverse systems
arising from~(\ref{ses:invlims}) yields the desired isomorphism.
\hfill
$\Box$
\vskip .2cm

So we can write $\FF^{\,(X,X_U)}$ as an inverse limit.  What have we
gained?  In the category of $\ZZ^n$-graded objects in which each graded
piece has finite dimension over $k$ (e.g.\ if the objects are chain
complexes which are finitely generated as $S$-modules), the functor
$^\Ast\invlim{}$ is exact, at least in the case where the inverse systems
are indexed by $\NN$---see \cite{Wei}, Exercise~3.5.2.  With this in mind
the following corollary is a simple consequence of \cite{Wei},
Theorem~3.5.8.

\begin{cor} \label{cor:homology}
To compute homology we have $H_i(\FF^{\,(X,X_U)}) =
{}^\Ast\invlim{t}H_i(\FF^{\,t}(X))$. \hfill $\Box$
\end{cor}

Until this point in this section, the labelled cell complex $X$ has been
arbitrary.  Now, however, we suppose that $X$ supports a cellular free
resolution of the ideal $I + \mm^{\aa + {\bf 1}}$, with $\aa \succeq
\aa_I$.  We will see shortly that for any $t$ the only nonvanishing
homology of $\FF^{\,t}(X)$ is in homological degree $1-n$, so the
previous corollary implies that the same holds for $\FF^{\,(X,X_U)}$.
Now $\FF_X$ has length at least $n-1$ (i.e. $\dim X \geq n-1$) because it
gives a free resolution of an artinian ideal; if we are so lucky that
$\FF_X$ has length exactly $n-1$, then the summand of $\FF^{\,(X,X_U)}$
in homological degree $1-n$ will be the last nonzero term.  In other
words, $\FF^{\,(X,X_U)}$ will be a \emph{free resolution} of some
$S$-module.  This is what makes Theorem~\ref{thm:relcocell} a special
case of the next result.  Even if we aren't so lucky with the length of
$\FF_X$, at least it will be split exact in homological degrees $> n-1$
(so that $\FF^{\,(X,X_U)}$ is split exact in homological degrees $<
1-n$), and we can still determine what the nonzero homology module is:
\begin{thm} \label{thm:coresolution}
Under the above conditions, $H_i(\FF^{\,(X,X_U)}) = 0$ whenever $\,i \neq
1-n$, and $H_{1-n}(\FF^{\,(X,X_U)}) = I^{[\aa]}[\aa + {\bf 1}]$.
\end{thm}
{\it Proof:\ } Let $J = I + \mm^{\aa + {\bf 1}}$.  For any $\bb \succeq
{\bf 0}$ Corollary~\ref{cor:fbX=SbI} implies that $\FF_\bb(X)$ is a free
resolution of the module $S[\bb]/\,S[\bb]\!\cap\!\widetilde{J}$,
homologically shifted down by 1.  Thus $\FF^\bb(X)$, which is the
$S$-dual of $\FF_\bb(X)$, is a complex whose homology in degree $i-1$ is
$\eext^i_S\Bigl( S[\bb]/ \,S[\bb]\!\cap\!\widetilde{J}\,,\,S\Bigr)$.  Now
$S[\bb]/\,S[\bb] \!\cap\!\widetilde{J} \subseteq T / \widetilde{J}$ is
artinian since $J = I + \mm^{\aa + {\bf 1}}$ is, and it is noetherian
because $S[\bb]$ is.  Hence the Ext module in question is, by \cite{BH},
Theorem~3.3.10(c), nonzero only for $i = n$.  Moreover,
Proposition~\ref{prop:ext} produces the equality
$$
  \eext^n_S \Bigl( S[\bb]/ \,S[\bb]\!\cap\!\widetilde{J}, S \Bigr) \ = \
  \Bigl(I\,/\:I\!\cap\!\mm^{\aa + \bb + {\bf 1}}\Bigr)[\aa + {\bf 1}].
$$
Taking the ${}^\Ast\invlim{\bb}$ of this last line and applying
Corollary~\ref{cor:homology} along with the completeness of $S$ proves
the theorem.
%

\end{section}

Department of Mathematics, University of California, Berkeley, CA
94720

e-mail: {\tt enmiller@math.berkeley.edu}

\end{document}